\documentclass[11pt]{amsart}
\usepackage{fullpage}
\usepackage{graphicx}
\usepackage{subfigure}
\usepackage{psfrag}
\usepackage{color}
\usepackage{amsmath}
\usepackage{amssymb}

\newtheorem{theorem}{Theorem}

\newtheorem{remark}[theorem]{Remark}

\title{A variational multiscale {N}ewton--{S}chur approach for the incompressible {N}avier--{S}tokes equations}

\author{D. Z. Turner} \author{K.~B.~Nakshatrala} \author{K.~D.~Hjelmstad}
\address{Correspondence to: Daniel Z. Turner, Department of Civil and Environmental Engineering, 
  2103 Newmark Laboratory, University of Illinois at Urbana-Champaign, Urbana, Illinois - 61801.}
\email{dzturne1@illinois.edu}
\address{Dr. Kalyana Babu Naskshatrala, Department of Civil and Environmental Engineering, 
  2524 Hydrosystems Laboratory, University of Illinois at Urbana-Champaign, Urbana, Illinois - 61801.}
\email{nakshatr@illinois.edu}
\address{Professor Keith D Hjelmstad, Department of Civil and Environmental Engineering, 
  3129e Newmark Laboratory, University of Illinois at Urbana-Champaign, Urbana, Illinois - 61801.}
\email{kdh@illinois.edu}

\date{\today}

\begin{document}

\begin{abstract}
In the following paper, we present a consistent Newton-Schur solution approach for variational multiscale formulations of the time-dependent Navier--Stokes equations in three dimensions.  The main contributions of this work are a systematic study of the variational multiscale method for three-dimensional problems, and an implementation of a consistent formulation suitable for large problems with high nonlinearity, unstructured meshes, and non-symmetric matrices. In addition to the quadratic convergence characteristics of a Newton-Raphson based scheme, the Newton-Schur approach increases computational efficiency and parallel scalability by implementing the tangent stiffness matrix in Schur's complement form.  As a result, more computations are performed at the element level.  Using a variational multiscale framework, we construct a two-level approach to stabilizing the incompressible Navier--Stokes equations based on a coarse and fine-scale subproblem. We then derive the Schur's complement form of the consistent tangent matrix. We demonstrate the performance of the method for a number of three-dimensional problems for Reynolds number up to 1000 including steady and time-dependent flows.  
\end{abstract}


\maketitle

\section{INTRODUCTION}

The incompressible Navier--Stokes equations are used to model a number of important physical phenomena including, pipe flow, flow around airfoils, ocean currents, weather, and blood flow in an artery, among others.  Significant emphasis has been placed in the literature on developing stabilized formulations robust enough to model complex flows at high Reynolds number \cite{Gunzburger, Gresho, Donea}. The variational multiscale method, due to Hughes \cite{Hughes2}, has been gaining popularity as a robust stabilization technique. This method consists of decomposing the problem into a coarse (or resolvable) and fine-scale (unresolvable) subproblem. The inclusion of the fine-scale influence stabilizes the solution.

The variational multiscale method has been well studied for two-dimensional problems, but less so for three-dimensional problems.  Also, three-dimensional effects become more influential as the Reynolds number increases.  It is well known that two-dimensional analysis cannot accurately capture flow features for high Reynolds flows. For example, it is well known that for the backward facing step problem at Reynolds number greater than 450, the reattachment length obtained by two-dimensional analysis does not match experimental results \cite{Jiang}. One of the main contributions of this paper is a systematic study of the variational multiscale method for three-dimensional problems.  To solve three-dimensional problems requires many more degrees of freedom. Using a mixed formulation, in which velocity and pressure are treated as unknowns, the number of degrees of freedom easily exceeds several million. Another contribution of this work is the Schur's complement implementation which greatly reduces the problem size and improves parallel performance.

The main challenge to numerically modeling the incompressible Navier--Stokes equations using the finite element method is due to the well known instability associated with using equal order interpolation for the velocity and pressure (which is computationally most convenient). Many stabilized methods for the incompressible Navier--Stokes equations have been developed that allow for equal order interpolation including \cite{Russo, Masud, Franca2, Hughes2, PSPG, SUPG}.  Stabilized formulations that employ a variational multiscale framework have been gaining popularity.  First introduced as a pathway to developing stabilized formulations from first principles, the variational multiscale concept has been used in a number of contexts, particularly for the incompressible Navier--Stokes equations \cite{TurnerStokes, TurnerNavierStokes, Nakshatrala, Nakshatrala2, Masud, Masud2, Koobus, Hughes3}. 

One can either solve the fine-scale problem in terms of the coarse-scale variables, thereby generating a stabilization term, and substitute back into the coarse-scale problem, or solve both problems simultaneously.  Typically a fixed-point iteration technique is used to linearize the weak form of the governing equations.  We have shown in \cite{TurnerNavierStokes} that by using a consistent Newton-Raphson based linearization, quadratic convergence may be achieved and that convergence may be achieved for problems for which fixed-point iterative methods diverge or converge very slowly. In this work, we present a Schur's complement implementation of the consistent Newton-Raphson based linearization of incompressible Navier--Stokes equations that provides for good computational efficiency and scalability.  We extend the work presented in \cite{TurnerNavierStokes} to large, nonlinear, three-dimensional problems on unstructured meshes.

Using the Schur's complement implementation of consistent formulations has been used primarily in the field of domain decomposition techniques and computational plasticity \cite{Kulkarni, Kulkarni2}. The Schur's complement implementation requires using block Gauss elimination to eliminate the fine-scale terms from the consistent tangent matrix.  The resulting tangent matrix contains modified coarse-scale terms that take into account the fine-scale influence.  If needed, the fine-scale terms can be obtained via post-processing. The Schur's complement form of the consistent tangent matrix reduces the problem size, and provides for much greater parallel scalability.

Primary emphasis in this paper is placed on solving several canonical three-dimensional problems for which experimental and numerical results are available in the literature. We show that using the Newton-Schur approach not only improves convergence, but also allows for much more efficient use of parallel resources. In the first section, we introduce the variational multiscale framework and describe the decomposition of the problem into a coarse and fine-scale subproblem.  In the next section, we present the consistent linearization of the Navier--Stokes equations and incorporate the Schur's complement form of the consistent tangent matrix to complete the Newton-Schur solution approach.  We then present a number of numerical examples that illustrate the advantages of the consistent Newton-Schur solution approach.
%
  \section{GOVERNING EQUATIONS}
  Let $\Omega$ be a bounded open domain, and $\Gamma$ be its boundary, 
  which is assumed to be piecewise smooth. Let the velocity vector field be denoted by 
  $\boldsymbol{v} : \Omega \rightarrow \mathbb{R}^{nd}$, where ``$nd$'' 
  is the number of spatial dimensions. Let the (kinematic) pressure field 
  be denoted by $p:\Omega \rightarrow \mathbb{R}$. 
  As usual, $\Gamma$ is divided into two parts, denoted by 
  $\Gamma^{\boldsymbol{v}}$ and $\Gamma^{\boldsymbol{t}}$, such 
  that $\Gamma^{\boldsymbol{v}} \cap \Gamma^{\boldsymbol{t}} =\emptyset$ 
  and $\Gamma^{\boldsymbol{v}} \cup \Gamma^{\boldsymbol{t}}=\Gamma$.  $\Gamma^{\boldsymbol{v}}$, and $\Gamma^{\boldsymbol{t}}$ are the Dirichlet and Neumann boundaries respectively. 
  The governing equations for incompressible Navier--Stokes flow can be written as 
  \begin{align}
    \label{Eqn:SC_Equilibrium}
    \boldsymbol{v} \cdot \nabla \boldsymbol{v} - 2\nu \nabla^2 \boldsymbol{v} + 
    \nabla p &= \boldsymbol{b} \qquad \; \;  \ \mbox{in}  \quad \Omega \\
    \label{Eqn:SC_Continuity}
    \nabla \cdot \boldsymbol{v} &= 0 \qquad \; \;   \ \mbox{in} \quad \Omega \\
    \label{Eqn:SC_VelocityBC}
    \boldsymbol{v} &= \boldsymbol{v}^{\mathrm{p}} 
    \qquad \; \mbox{on} \quad \Gamma^{\boldsymbol{v}} \\
    \label{Eqn:SC_TractionBC}
    -p \boldsymbol{n} + \nu (\boldsymbol{n} \cdot \nabla) \boldsymbol{v} &= 
    \boldsymbol{t}^{\boldsymbol{n}} \qquad \; \mbox{on} \quad \Gamma^{\boldsymbol{t}} 
  \end{align}
  where $\boldsymbol{v}$ is the velocity, $p$ is the kinematic pressure (pressure divided by density),  $\nabla$ is the gradient operator, $\nabla^2$ is the Laplacian operator, $\boldsymbol{b}$ is the body force, $\nu > 0$ is the kinematic viscosity, $\boldsymbol{v}^{\mathrm{p}}$ is the 
  prescribed velocity vector field, $\boldsymbol{t}^{\boldsymbol{n}}$ is the 
  prescribed traction, and $\boldsymbol{n}$ is the unit outward normal vector 
  to $\Gamma$. Equation \eqref{Eqn:SC_Equilibrium} represents the balance of 
  linear momentum, and equation \eqref{Eqn:SC_Continuity} represents the continuity 
  equation for an incompressible continuum. Equations \eqref{Eqn:SC_VelocityBC} 
  and \eqref{Eqn:SC_TractionBC} are the Dirichlet and Neumann boundary 
  conditions, respectively. 
 
%
%
\section{VARIATIONAL MULTISCALE FRAMEWORK}
It is well known that the classical mixed formulation for the incompressible Navier--Stokes equations does not produce stable numerical results. The instability of the standard Galerkin formulation may be mathematically explained by the LBB \emph{inf-sup} stability condition \cite{Brezzi, Gunzburger}. To get stable numerical results, one must stabilize the standard Galerkin formulation.  The variational multiscale concept, which stems from the pioneering work by Hughes \cite{Hughes2}, decomposes the underlying fields into coarse or resolvable scales and subgrid or unresolvable scales.
%
\subsection{Multiscale decomposition}
Let us divide the domain $\Omega$ into $N$ non-overlapping subdomains 
$\Omega^e$ (which in the finite element context will be elements) such 
that
%
\begin{equation}
  \Omega = \overset{N}{\underset{e = 1}{\bigcup}} \Omega^{e}
\end{equation}
The boundary of the element $\Omega^{e}$ is denoted by $\Gamma^{e}$. We 
decompose the velocity field $\boldsymbol{v}(\boldsymbol{x})$ into 
coarse-scale and fine-scale components, indicated as $\bar{\boldsymbol{v}}
(\boldsymbol{x})$ and $\boldsymbol{v}'(\boldsymbol{x})$, respectively.  
%
\begin{equation}
  \label{Eqn:SC_Decompose_v}
  \boldsymbol{v}(\boldsymbol{x}) = \bar{\boldsymbol{v}} 
  (\boldsymbol{x}) + \boldsymbol{v}'(\boldsymbol{x})
\end{equation}
Likewise, we decompose the weighting function $\boldsymbol{w}(\boldsymbol{x})$ 
into coarse-scale $\bar{\boldsymbol{w}}(\boldsymbol{x})$ and fine-scale 
$\boldsymbol{w}'(\boldsymbol{x})$ components.
%
\begin{equation}
  \label{Eqn:SC_Decompose_w}
  \boldsymbol{w}(\boldsymbol{x}) = \bar{\boldsymbol{w}}(\boldsymbol{x}) 
  + \boldsymbol{w}'(\boldsymbol{x})
\end{equation}
We further make an assumption that the fine-scale 
components vanish along each element boundary.
%
\begin{equation}
  \label{Eqn:SC_Fine_Scale_Vanish}
  \boldsymbol{v}'(\boldsymbol{x}) = \boldsymbol{w}'(\boldsymbol{x}) 
  = \boldsymbol{0} \quad \mbox{on} \quad \Gamma^{e} \; ; \; e = 1, \ldots, N
\end{equation}
Let $\mathcal{V}$ be the function space for the coarse-scale component 
of the velocity $\bar{\boldsymbol{v}}$, and $\mathcal{W}$ be the 
function space for $\bar{\boldsymbol{w}}$; and are defined as 
\begin{align}
  \label{Eqn:SC_Function_Space_for_v}
  \mathcal{V} &:= \{\boldsymbol{v}  \; \big| \; \boldsymbol{v} \in (H^{1}(\Omega))^{nd},
  \boldsymbol{v} = \boldsymbol{v}^{\mathrm{p}} \; \mathrm{on} \; \Gamma^{\boldsymbol{v}} \} \\
  \label{Eqn:SC_Function_Space_for_w}
  \mathcal{W} &:= \{\boldsymbol{w}  \; \big| \; \boldsymbol{w} \in (H^{1}(\Omega))^{nd},
  \boldsymbol{w} = \boldsymbol{0} \; \mathrm{on} \; \Gamma^{\boldsymbol{v}} \}
\end{align} 
Let $\mathcal{V}'$ be the function space for both the fine-scale 
component of the velocity $\boldsymbol{v}'$ and its corresponding weighting 
function $\boldsymbol{w}'$, and is defined as 
\begin{align}
  \label{Eqn:SC_function_space_for_v_prime}
  \mathcal{V}' := \{\boldsymbol{v} \; \bigr| \; \boldsymbol{v} 
  \in (H^{1}(\Omega^e))^{nd} , \; \boldsymbol{v} = \boldsymbol{0} \; \mbox{on} \; 
  \Gamma^e, e = 1, \ldots, N \} 
\end{align}
In theory, we could decompose the pressure field into coarse-scale and fine-scale 
components. However, for simplicity we assume that there are no fine-scale terms 
for the pressure $p(\boldsymbol{x})$ and for its corresponding weighting function 
$q(\boldsymbol{x})$. Hence, the function space for the fields $p(\boldsymbol{x})$ 
and $q(\boldsymbol{x})$ is $\mathcal{P}$
\begin{equation}
  \label{Eqn:SC_Function_Space_for_p}
  \mathcal{P} := \{p \; \big| \; p \in L^{2}(\Omega) \}
\end{equation}
where $L^2(\Omega)$ is the space of square-integrable functions on the domain $\Omega$. For further details on function spaces refer to Brezzi and Fortin \cite{Brezzi}.
%
\subsection{Two-level classical mixed formulation}
After substitution of equations \eqref{Eqn:SC_Decompose_v} and \eqref{Eqn:SC_Decompose_w} 
into the classical mixed formulation and decomposing the problem into a coarse and fine-scale subproblem, the resulting weak form can be written as two parts. The \emph{coarse-scale problem} can be written as:
\begin{align}
  \label{Eqn:SC_CoarseScale_Momentum}
(\bar{\boldsymbol{w}}, (\bar{\boldsymbol{v}} + \boldsymbol{v}') \cdot \nabla (\bar{\boldsymbol{v}} + \boldsymbol{v}')) + (\nabla \bar{\boldsymbol{w}},2\nu \nabla (\bar{\boldsymbol{v}}+\boldsymbol{v}')) - (\nabla \cdot \bar{\boldsymbol{w}}, p)  &= (\bar{\boldsymbol{w}} , \boldsymbol{b} ) + ( \bar{\boldsymbol{w}}, \boldsymbol{h})_{\Gamma_t} \quad \forall \ \bar{\boldsymbol{w}} \in \bar{\mathcal{W}} \\
  \label{Eqn:SC_CoarseScale_Continuity}
  (q, \nabla \cdot (\bar{\boldsymbol{v}} + \boldsymbol{v}')) &= 0  \quad \forall \ q \in \mathcal{P}
\end{align}
The \emph{fine-scale problem} can be written as:
\begin{align}
  \label{Eqn:SC_FineScale_Momentum}
  (\boldsymbol{w}', (\bar{\boldsymbol{v}} + \boldsymbol{v}') \cdot \nabla (\bar{\boldsymbol{v}} + \boldsymbol{v}')) + (\nabla  \boldsymbol{w}',2\nu \nabla (\bar{\boldsymbol{v}}+\boldsymbol{v}')) - (\nabla \cdot  \boldsymbol{w}', p) &= (\boldsymbol{w}', \boldsymbol{b} ) + (\boldsymbol{w}', \boldsymbol{h})_{\Gamma_t} \quad \forall \ \boldsymbol{w}' \in \mathcal{W}' 
\end{align}
For convenience, we define the $L^2$ inner-product over a spatial domain, $K$, as
\begin{align}
  (\boldsymbol{a}, \boldsymbol{b})_K = \int_K \boldsymbol{a} \cdot \boldsymbol{b} \; \mathrm{d} K
\end{align}
The subscript $K$ will be dropped if $K$ is the whole of $\Omega$, that is $K = \Omega$. Notice that the weak form is nothing more than the classical mixed formulation written for the coarse and fine-scale varaibles, with the exception that there is no weak form of the fine-scale incompressibility constraint. For a more detailed derivation of the weak form, see \cite{TurnerNavierStokes}. 

%
\subsection{Fine--scale interpolation and bubble functions}
If one chooses a single bubble function for interpolating the fine-scale variables (similar to the MINI element \cite{Arnold}), then we have
\begin{align}
 \label{eq:DefineBubbles3}
\boldsymbol{v}' = b^e \boldsymbol{\beta}; \quad \boldsymbol{w}' = b^e \boldsymbol{\gamma}
\end{align}
where $b^e$ is a bubble function, and $\boldsymbol{\beta}$ and $\boldsymbol{\gamma}$ 
are constant vectors. 
The gradients of the fine-scale velocity and weighting functions are
\begin{align}
\nabla \boldsymbol{v}' = \boldsymbol{\beta}\nabla b^{eT}; \quad \nabla \boldsymbol{w}' = \boldsymbol{\gamma}\nabla b^{eT}
\end{align}
where $\nabla b^e $ is a $\mathrm{dim} \times 1 $ vector of the derivatives of the bubble function.
%
%
%

%
\section{CONSISTENT NEWTON--SCHUR SOLUTION STRATEGY}
For the consistent Newton--Schur method, we treat equations \eqref{Eqn:SC_CoarseScale_Momentum}--\eqref{Eqn:SC_FineScale_Momentum} as global residuals and obtain the solution using an iterative update equation.
\subsection{Vector residual}
First, we express the solution field and its weighting function in terms of nodal values $\bar{\boldsymbol{v}} = \hat{\bar{\boldsymbol{v}}}^T \boldsymbol{N}^T$ and $\bar{\boldsymbol{w}} = \hat{\bar{\boldsymbol{w}}}^T \boldsymbol{N}^T$, respectively (where $\boldsymbol{N}$ is a row vector of shape functions for each node). After substituting these expressions into the global residuals, the vector residuals, $\boldsymbol{R}$, that are the sum contributions of the vector residuals at the element level, $\boldsymbol{R}^e$, can be written as
\begin{align}
  \label{Eqn:SC_Coarse_Resid_Vec}
  \boldsymbol{R}^e_c(\bar{\boldsymbol{v}};\boldsymbol{v}',p)  := &\int_{\Omega} (\boldsymbol{N}^T \odot \boldsymbol{I}) ((\bar{\boldsymbol{v}} + \boldsymbol{v}') \cdot \nabla) 
(\bar{\boldsymbol{v}} + \boldsymbol{v}') \; \mathrm{d} \Omega  \nonumber \\
&+ \;2\nu \int_{\Omega}  ((\boldsymbol{DN})\boldsymbol{J}^{-1} \odot \boldsymbol{I})
\mathrm{vec}[\nabla \bar{\boldsymbol{v}} + \nabla \boldsymbol{v}'] \; \mathrm{d} \Omega \nonumber \\
& - \int_{\Omega} \mathrm{vec}[\boldsymbol{J}^{-T}(\boldsymbol{DN})^T]p \; \mathrm{d} \Omega - \int_{\Omega} (\boldsymbol{N}^T \odot \boldsymbol{I}) \boldsymbol{b} \; \mathrm{d} \Omega  - \int_{\Gamma_t} (\boldsymbol{N}^T \odot \boldsymbol{I}) \boldsymbol{h} \; \mathrm{d} \Gamma\\
  \label{Eqn:SC_Coarse_Resid_Vec_2}
  \boldsymbol{R}^e_p(\bar{\boldsymbol{v}};\boldsymbol{v}')  := & - \ \int_{\Omega} \boldsymbol{N}^T \nabla \cdot (\bar{\boldsymbol{v}} + \boldsymbol{v}') \; \mathrm{d} \Omega  \\
  \label{Eqn:SC_Fine_Resid_Vec}
  \boldsymbol{R}^e_f(\bar{\boldsymbol{v}};\boldsymbol{v}',p)  := & \int_{\Omega_e} b^e \boldsymbol{I}((\bar{\boldsymbol{v}} + \boldsymbol{v}') \cdot \nabla) 
(\bar{\boldsymbol{v}} + \boldsymbol{v}') \; \mathrm{d} \Omega + \;2\nu \int_{\Omega_e}  (\nabla b^{eT} \odot \boldsymbol{I})
\mathrm{vec}[\nabla \bar{\boldsymbol{v}} + \nabla \boldsymbol{v}'] \; \mathrm{d} \Omega \nonumber \\
&  - \int_{\Omega_e}  (\nabla b^{eT} \odot \boldsymbol{I})\mathrm{vec}[\boldsymbol{I}]p \; \mathrm{d} \Omega - \int_{\Omega_e}  b^e \boldsymbol{I} \boldsymbol{b} \; \mathrm{d} \Omega 
\end{align}
where the subscripts `$c$', `$p$', and `$f$' stand for \emph{coarse}, \emph{pressure}, and \emph{fine}, and $\boldsymbol{DN}$ represents a matrix of the first derivatives of the element shape functions, which is defined for a triangular element as
\begin{align}
 \boldsymbol{DN} &:= \left[ \begin{array}{cc} \frac{\partial N_1}{\partial \xi_1} & \frac{\partial N_1}{\partial \xi_2} \\
\vdots & \vdots \\
\frac{\partial N_3}{\partial \xi_1} & \frac{\partial N_3}{\partial \xi_2} \end{array} \right]
\end{align}
$\boldsymbol{J}$ is the element jacobian matrix, $\mathrm{vec}[\cdot]$ is an operation that represents a matrix with a vector, $\odot$ is the Kronecker product \cite{Graham}. 

\begin{remark}
For transient problems, we simply add the time derivative of the velocity to equation \eqref {Eqn:SC_Equilibrium}.  Using the backward Euler method to integrate in time from step $n$ to step $n+1$, we add the following terms, $\int_{\Omega} \frac{1}{\Delta t}(\boldsymbol{N}^T \odot \boldsymbol{I}) \bar{\boldsymbol{v}}_{n+1} \; \mathrm{d} \Omega 
  - \int_{\Omega} \frac{1}{\Delta t}(\boldsymbol{N}^T \odot \boldsymbol{I}) \bar{\boldsymbol{v}}_{n} \; \mathrm{d} \Omega$ and $\int_{\Omega_e} \frac{b^e}{\Delta t} \boldsymbol{I} \bar{\boldsymbol{v}}_{n+1} \; \mathrm{d} \Omega
  - \int_{\Omega_e} \frac{b^e}{\Delta t} \boldsymbol{I} \bar{\boldsymbol{v}}_{n} \; \mathrm{d} \Omega$ to equations \eqref{Eqn:SC_Coarse_Resid_Vec} and \eqref{Eqn:SC_Fine_Resid_Vec}, respectively, where $\bar{\boldsymbol{v}}_n$ is the converged velocity from the previous time step. For furhter details, see \cite{TurnerNavierStokes}.
\end{remark}
\subsection{Tangent matrix}
Using a Newton-Raphson type approach, we obtain the solution in an iterative fashion using the following update equation until the residual is under a prescribed tolerance,
\begin{align}
\bar{\boldsymbol{v}}^{i+1} = \bar{\boldsymbol{v}}^i + \Delta \bar{\boldsymbol{v}}^i \quad ; \quad p^{i+1} = p^i + \Delta p^i \quad ; \quad \boldsymbol{v}'^{i+1} = \boldsymbol{v}'^i + \Delta \boldsymbol{v}'^i
\end{align}
where the updates at each iteration, $i$ are calculated from the following system of equations
\begin{align}
\label{eq:NewtonUpdate2}
\left[ \begin{array}{ccc}   \frac{\textrm{D}\boldsymbol{R}_c}{\textrm{D}\bar{\boldsymbol{v}}} &\frac{\textrm{D}\boldsymbol{R}_c}{\textrm{D}p} & \frac{\textrm{D}\boldsymbol{R}_c }{\textrm{D}\boldsymbol{v}'} \\
  \frac{\textrm{D}\boldsymbol{R}_p}{\textrm{D}\bar{\boldsymbol{v}}} &  \frac{\textrm{D}\boldsymbol{R}_p}{\textrm{D}p} & \frac{\textrm{D}\boldsymbol{R}_p }{\textrm{D}\boldsymbol{v}'}  \\
 \frac{\textrm{D}\boldsymbol{R}_f}{\textrm{D}\bar{\boldsymbol{v}}} & \frac{\textrm{D}\boldsymbol{R}_f}{\textrm{D}p} & \frac{\textrm{D}\boldsymbol{R}_f}{\textrm{D}\boldsymbol{v}'} \end{array} \right] \left\{ \begin{array}{c} \Delta \bar{\boldsymbol{v}} \\
  \Delta p \\
  \Delta \boldsymbol{v}' \\ \end{array} \right\} = - \left\{ \begin{array}{c}  \boldsymbol{R}_c \\
   \boldsymbol{R}_p\\
   \boldsymbol{R}_f \\ \end{array} \right\}
\end{align}
The element matrices, $\frac{\textrm{D}\boldsymbol{R}^e}{\textrm{D}(\cdot)}$, which are assembled to form the consistent tangent matrix, are presented in \cite{TurnerNavierStokes}.  
%
\section{SCHUR'S COMPLEMENT IMPLEMENTATION}
 
By applying block Gauss elimination to equation \eqref{eq:NewtonUpdate2}, we obtain the Schur complement representation of the update equation,  
\begin{align}
\label{eq:SchurUpdate}
  \left[ \begin{array}{cc}  \boldsymbol{K}_{\bar{\boldsymbol{v}}\bar{\boldsymbol{v}}} & \boldsymbol{K}_{\bar{\boldsymbol{v}}p} \\
   \boldsymbol{K}_{p\bar{\boldsymbol{v}}} & \boldsymbol{K}_{pp} \end{array} \right] \left\{ \begin{array}{c} \Delta \bar{\boldsymbol{v}} \\
  \Delta p \\ \end{array} \right\} =  -\left\{ \begin{array}{c}   \boldsymbol{R}_1 \\
  \boldsymbol{R}_2  \end{array} \right\}
\end{align}
where each of the stiffness terms, $\boldsymbol{K}$, represent the assembled element contributions as follows:
\begin{align}
 \boldsymbol{K}_{\bar{\boldsymbol{v}}\bar{\boldsymbol{v}}}  = & \mathop{{\huge\boldsymbol{\mathsf{A}}}}_{e=1}^{Nele} \left[ \frac{\textrm{D}\boldsymbol{R}_{c}^e}{\textrm{D}\bar{\boldsymbol{v}}} - \frac{\textrm{D}\boldsymbol{R}_{c}^e }{\textrm{D}\boldsymbol{v}'}\left[\frac{\textrm{D}\boldsymbol{R}_f^e}{\textrm{D}\boldsymbol{v}'}\right]^{-1}\frac{\textrm{D}\boldsymbol{R}_f^e}{\textrm{D}\bar{\boldsymbol{v}}}\right] \nonumber \\
 \boldsymbol{K}_{\bar{\boldsymbol{v}}p}  = & \mathop{{\huge\boldsymbol{\mathsf{A}}}}_{e=1}^{Nele} \left[\frac{\textrm{D}\boldsymbol{R}_{c}^e}{\textrm{D}p} - \frac{\textrm{D}\boldsymbol{R}_{c}^e }{\textrm{D}\boldsymbol{v}'}\left[\frac{\textrm{D}\boldsymbol{R}_f^e}{\textrm{D}\boldsymbol{v}'}\right]^{-1}\frac{\textrm{D}\boldsymbol{R}_f^e}{\textrm{D}p} \right]\nonumber \\
 \boldsymbol{K}_{p\bar{\boldsymbol{v}}}  = & \mathop{{\huge\boldsymbol{\mathsf{A}}}}_{e=1}^{Nele} \left[\frac{\textrm{D}\boldsymbol{R}_{p}^e}{\textrm{D}\bar{\boldsymbol{v}}} - \frac{\textrm{D}\boldsymbol{R}_{p}^e }{\textrm{D}\boldsymbol{v}'}\left[\frac{\textrm{D}\boldsymbol{R}_f^e}{\textrm{D}\boldsymbol{v}'}\right]^{-1}\frac{\textrm{D}\boldsymbol{R}_f^e}{\textrm{D}\bar{\boldsymbol{v}}} \right]\nonumber \\
 \boldsymbol{K}_{pp}  = &  \mathop{{\huge\boldsymbol{\mathsf{A}}}}_{e=1}^{Nele} \left[- \frac{\textrm{D}\boldsymbol{R}_{p}^e}{\textrm{D}\boldsymbol{v}'}\left[\frac{\textrm{D}\boldsymbol{R}_f^e}{\textrm{D}\boldsymbol{v}'}\right]^{-1}\frac{\textrm{D}\boldsymbol{R}_f^e}{\textrm{D}p} \right]
\end{align}
And the residuals are given as
\begin{align}
 \boldsymbol{R}_1 & =  \mathop{{\huge\boldsymbol{\mathsf{A}}}}_{e=1}^{Nele} \left[\boldsymbol{R}_{c}^e - \frac{\textrm{D}\boldsymbol{R}_{c}^e }{\textrm{D}\boldsymbol{v}'}\left[\frac{\textrm{D}\boldsymbol{R}_f^e}{\textrm{D}\boldsymbol{v}'}\right]^{-1}\boldsymbol{R}_f^e \right]\nonumber \\
 \boldsymbol{R}_2 & =  \mathop{{\huge\boldsymbol{\mathsf{A}}}}_{e=1}^{Nele} \left[\boldsymbol{R}_{p}^e - \frac{\textrm{D}\boldsymbol{R}_{p}^e }{\textrm{D}\boldsymbol{v}'}\left[\frac{\textrm{D}\boldsymbol{R}_f^e}{\textrm{D}\boldsymbol{v}'}\right]^{-1}\boldsymbol{R}_f^e \right]
\end{align}
From equation \ref{eq:SchurUpdate} we update the coarse-scale velocity and pressure. We then solve for the fine-scale increment from
\begin{align}
 \Delta \boldsymbol{v}' = \left[\frac{\textrm{D}\boldsymbol{R}_f}{\textrm{D}\boldsymbol{v}'}\right]^{-1} \left[ -\boldsymbol{R}_f - \frac{\textrm{D}\boldsymbol{R}_f}{\textrm{D}\bar{\boldsymbol{v}}}\Delta \bar{\boldsymbol{v}} - \frac{\textrm{D}\boldsymbol{R}_f}{\textrm{D}p}\Delta p \right]
\end{align}
Again, the coarse and fine-scale velocity and pressure are updated until convergence.
In addition to faster convergence properties obtained for the Newton-Raphson type solution procedure, the Schur complement representation increases parallel performance by taking advantage of the fine-scale terms vanishing on the boundary of the element.  The fine-scale computations can therefore take place at the element level without the need for inter-processor communication.
%
\section{NUMERICAL RESULTS}

The above formulation was implemented in C++ on the Turing cluster which consists of 768 Apple Xservers, each with two 2 GHz G5 processors and 4 GB of RAM, for a total of 1536 processors \cite{Turing}. The cluster is connected using a high-bandwidth, low-latency Myrinet network from Myricom. The operating system used is Mac OS X Server, version 10.3. The data structures were implemented in parallel using the Portable, Extensible Toolkit for Scientific Computation (PETSc) \cite{petsc-web-page} Vec and Mat objects (note that the Vec object is different than the $\mathrm{vec}[\cdot]$ operator). PETSc, in turn uses MPICH for parallel communication. Linear solutions for each iteration were obtained using the iterative Krylov subspace solver, KSP, provided by PETSc. The tolerance was set to 1E-12 and the maximum number of iterations was set to 50. The code aborts if the max number of iterations are reached without meeting the tolerance.

Between one and 128 processors were used for each simulation and the results were visualized using Tecplot 360 \cite{Tecplot} and VisIt 1.8.1 \cite{Visit}. The mesh was partitioned using a simple block partitioning strategy. The parallel results may be improved with the use of a more sophisticated partitioning algorithm, but that is beyond the scope of this work. In order to avoid deadlock, matrix assembly routines are called by each process after all of its elements have been computed. 

\subsection{Body force driven cavity}
\begin{figure}[htb!]
	\centering
  \includegraphics[scale=0.5]{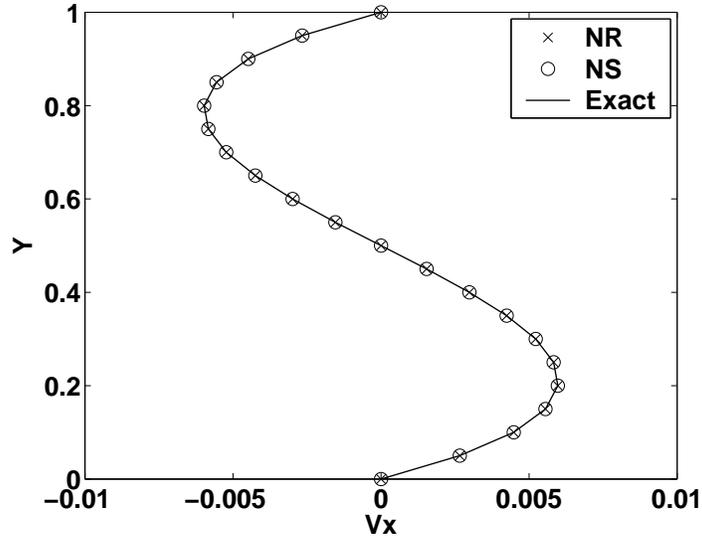}
	\caption{Body force driven cavity: velocity in the $x$-direction along the centerline of the cavity for the Newton-Schur (NS) and Newton-Raphson (NR) solution approaches.}
	\label{fig:2DExactSolVel}
\end{figure}

Using the body force driven cavity problem, we compare the speedup on a single processor between a consistent formulation using the Schur's complement implementation and one without. The problem description, geometry, and boundary conditions for the body force driven cavity problem are given in \cite{TurnerNavierStokes}. An exact solution is derived for a given body force. The results for both the Newton-Raphson approach, presented in \cite{TurnerNavierStokes}, and the Newton--Schur approach, presented here, are shown in Figure \ref{fig:2DExactSolVel}. Figure \ref{fig:2DExactSolVel} shows the velocity in the $x$-direction along the center of the cavity as y increases from the bottom to the top. Notice that the results are identical for both methods.

\begin{figure}[htb!]
	\centering
  \includegraphics[scale=0.5]{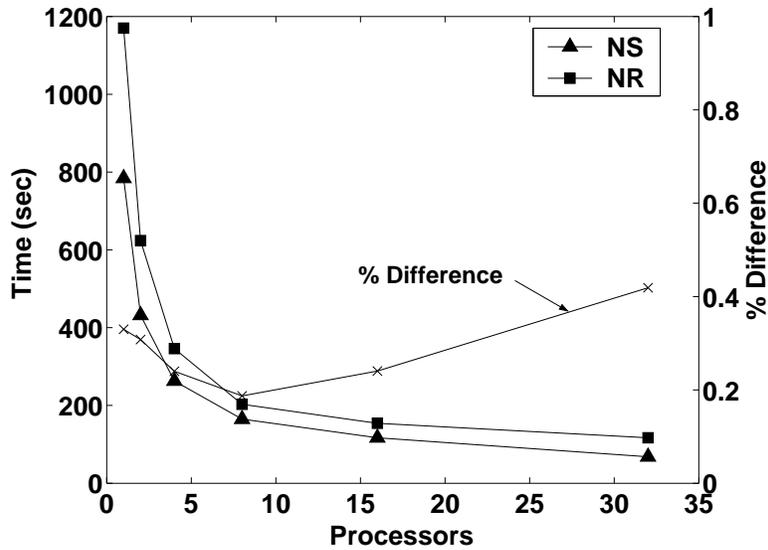}
	\caption{Body force driven cavity: computational time in seconds comparing the Newton-Schur (NS) and Newton-Raphson (NR) solution approaches for a given number of processors and a fixed problem size of 200,000 degrees of freedom.}
	\label{fig:CompTimeJet}
\end{figure}

Figure \ref{fig:CompTimeJet} compares the computational time required by the Newton--Raphson and Newton--Schur approaches as the number of processors increases. Notice that for 32 processors, the computational costs are reduced by over 40\%. Since the computational advantages of the Schur's complement implementation are based largely on parallel features, the benefits of the Newton--Schur approach grow as the number of processors increases.

\subsection{Three-dimensional lid-driven cavity}
\begin{figure}[htb!]
	\centering
  \includegraphics[scale=0.5]{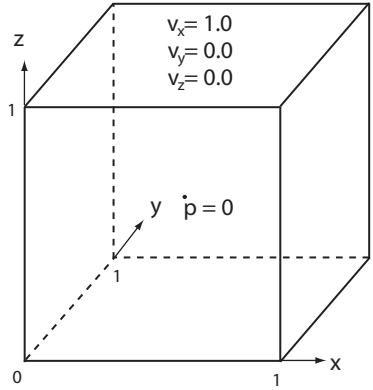}
	\caption{3D lid-driven cavity: problem description and boundary conditions. All boundaries are given no-slip wall conditions ($v_x = v_y = v_z = 0$) except for the top boundary which is given a unit velocity in the $x$-direction and zero velocity in $y$ and $z$. The pressure at one node at the center of the domain is prescribed as zero.}
	\label{fig:3DCavity}
\end{figure}

The geometry and boundary conditions for the three-dimensional lid-driven cavity are shown in Figure \ref{fig:3DCavity}. A unit velocity is prescribed in the positive $x$-direction across the top of a cube of unit volume. The viscosity is adjusted to obtain higher Reynolds numbers. Steady solutions to the three-dimensional lid-driven cavity have only been shown to exist for Reynolds number less than 2000. Experimental studies of the three-dimensional lid-driven cavity show a number of unsteady motions for $Re > 2000$, for example eddies and G$\ddot{\mathrm{o}}$rtler vortices \cite{Koseff}.

\begin{figure}[htb!]
	\centering
  \includegraphics[scale=0.5]{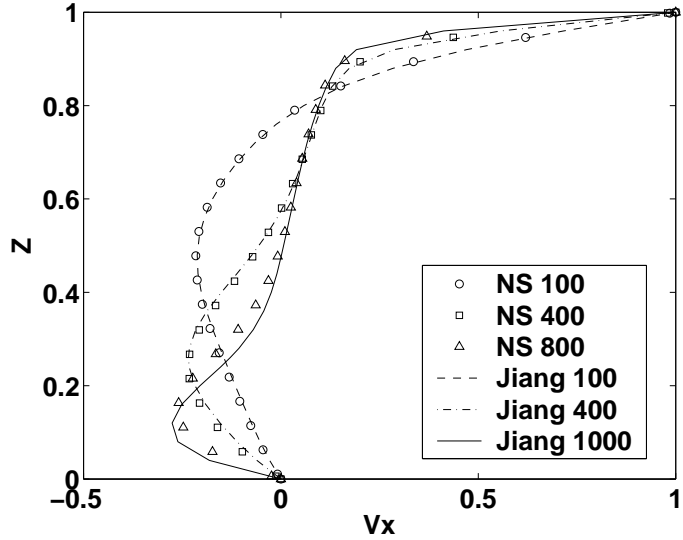}
	\caption{3D lid-driven cavity: velocity in the $x$-direction as measured along the center of the cavity in the $z$-direction ($x=0,y=0,z$) compared to the results obtained in \cite{Jiang}.}
	\label{fig:CenterlineCavity3D}
\end{figure}

Figure \ref{fig:CenterlineCavity3D} shows the velocity in the $x$-direction along the centerline of the cavity along the height of the cavity. Notice that the results correspond well with those reported in \cite{Jiang}. Using the Newton--Schur solution approach we were able to obtain steady state solutions up to $Re = 865$.  Above $Re = 865$ the Newton--Schur solution approach diverges. In \cite{Jiang} steady-state solutions are obtained for $Re = 1000$ using a velocity-pressure-vorticity formulation and a time-marching scheme. To our knowledge, this paper presents the first results for high Reynolds flows for the lid-driven cavity without the use of a time-marching scheme. As the Reynolds number increases, the instabilities created by boundary layer flow over concave surfaces, like the corners of the cavity, begin to take effect. These instabilities represent the transition into turbulent flow \cite{Saric}. The failure to converge above $Re = 800$ for the Newton--Schur solution strategy could be due to the presence of Taylor-G$\ddot{\mathrm{o}}$rtler vortices as reported in \cite{Jiang}.

Taylor-G$\ddot{\mathrm{o}}$rtler vortices occur when the velocity profile approaches zero at the boundary over a concave (and in some cases convex) surface. They are the result of a centrifugal instability. Rayleigh first recognized this instability in 1916 and showed that as the distance increases radially from the center point of a curved surface, the velocity must also increase, otherwise an inviscid axisymmetric instability occurs. In the case of the corners of the three-dimensional cavity, the flow decreases due to friction with the walls causing a decrease in velocity with radial distance. 

\begin{figure}[htb!]
	\centering
  \subfigure{\includegraphics[scale=0.2]{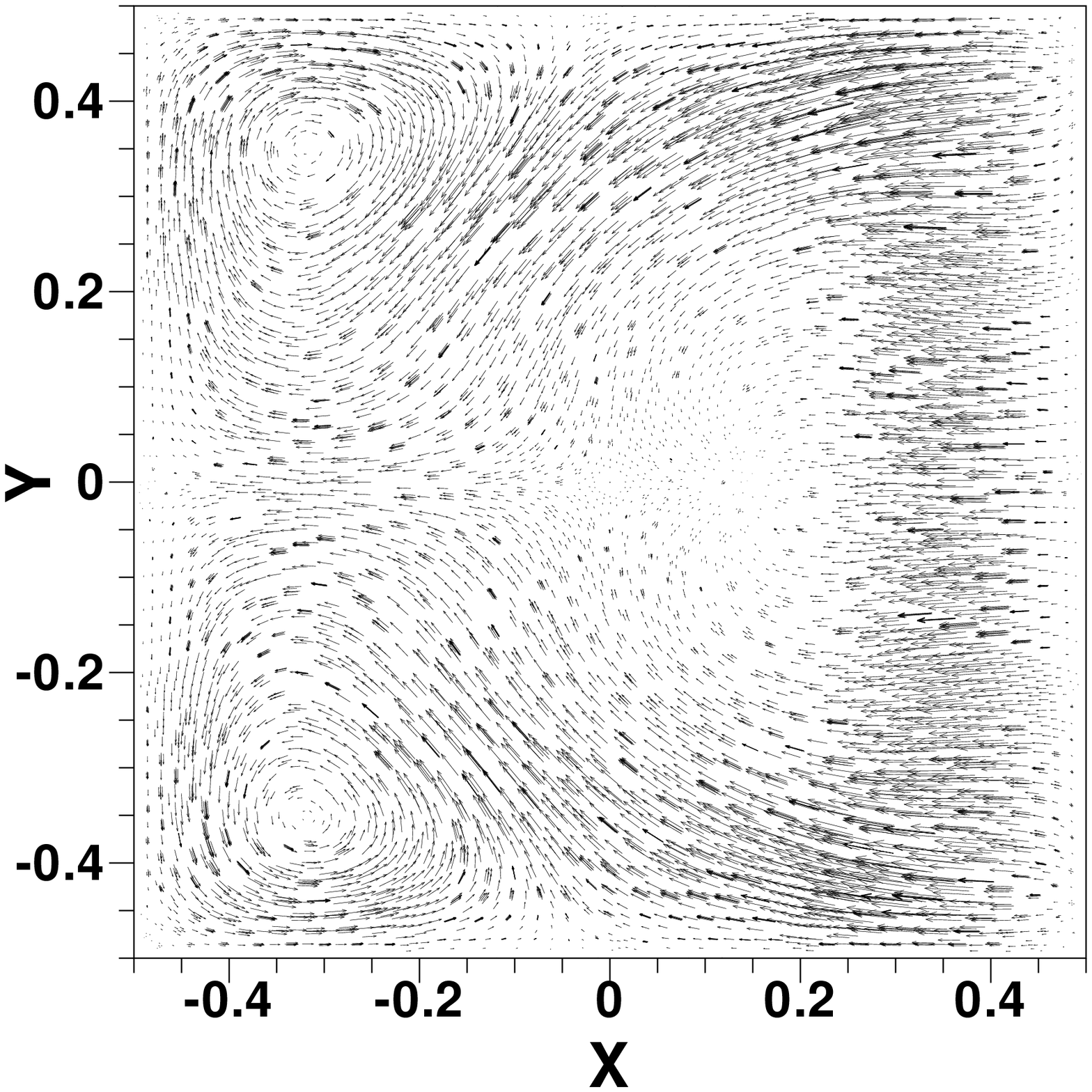}}
  \subfigure{\includegraphics[scale=0.2]{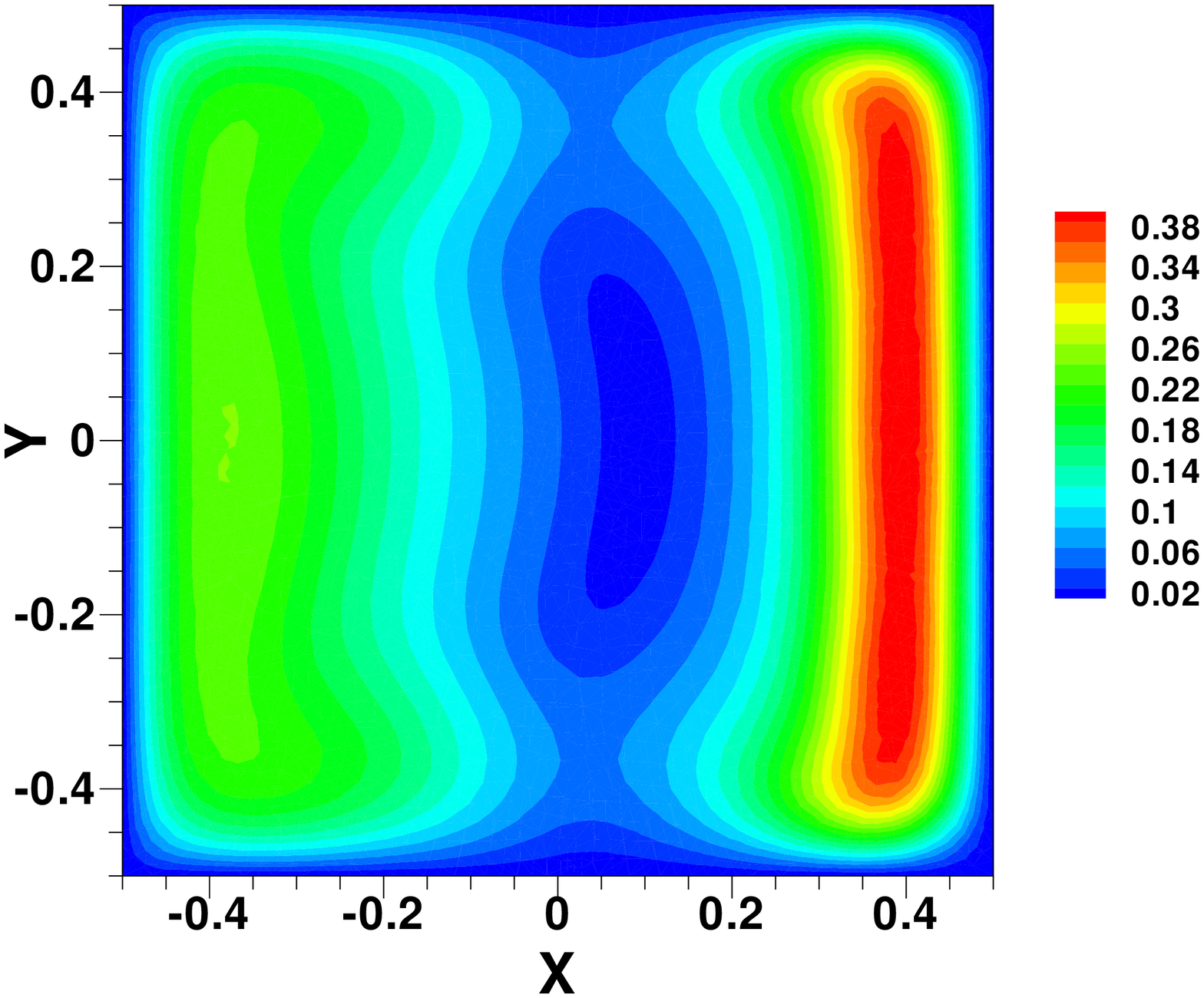}}
  \subfigure{\includegraphics[scale=0.2]{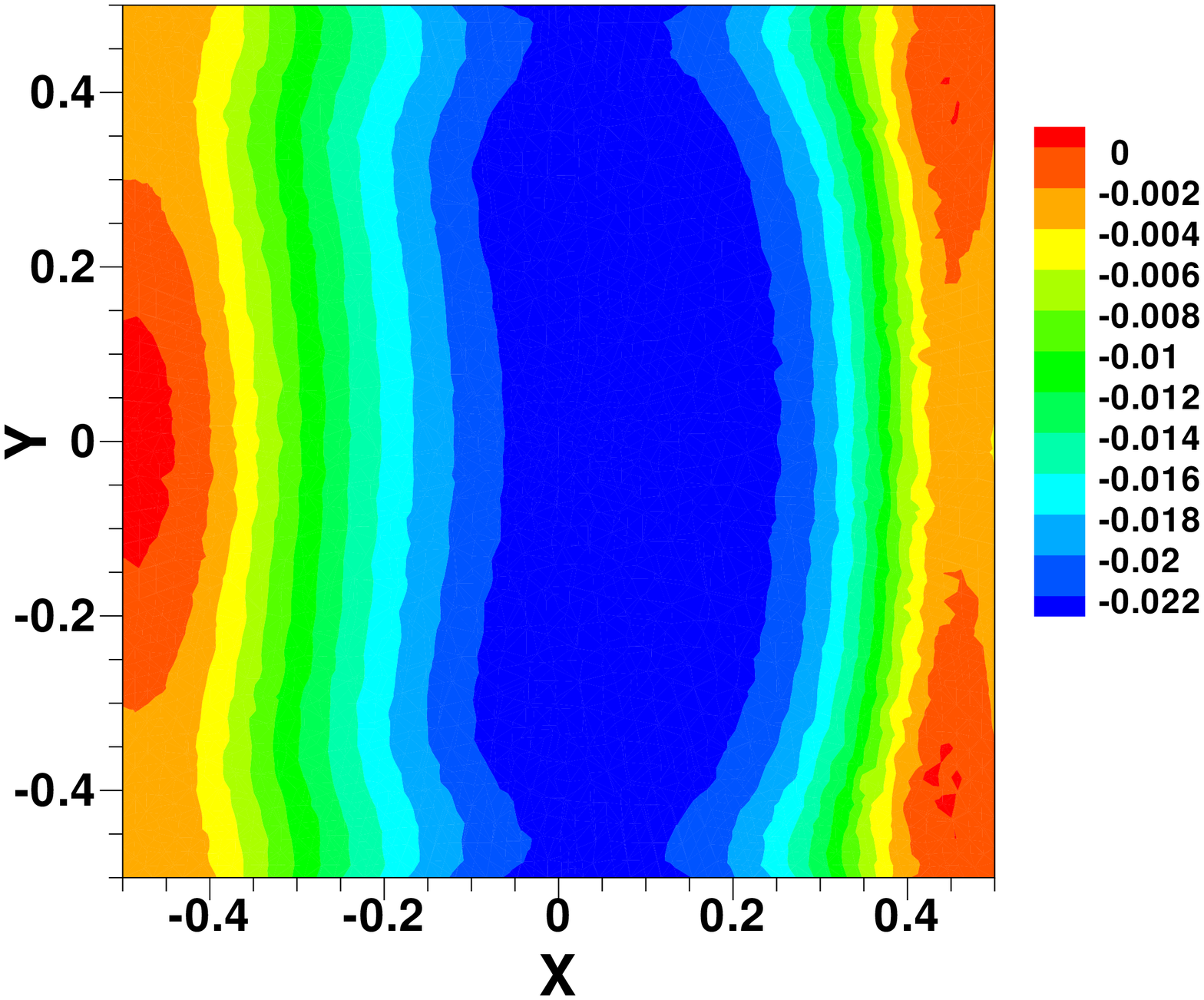}}  
    \subfigure{\includegraphics[scale=0.2]{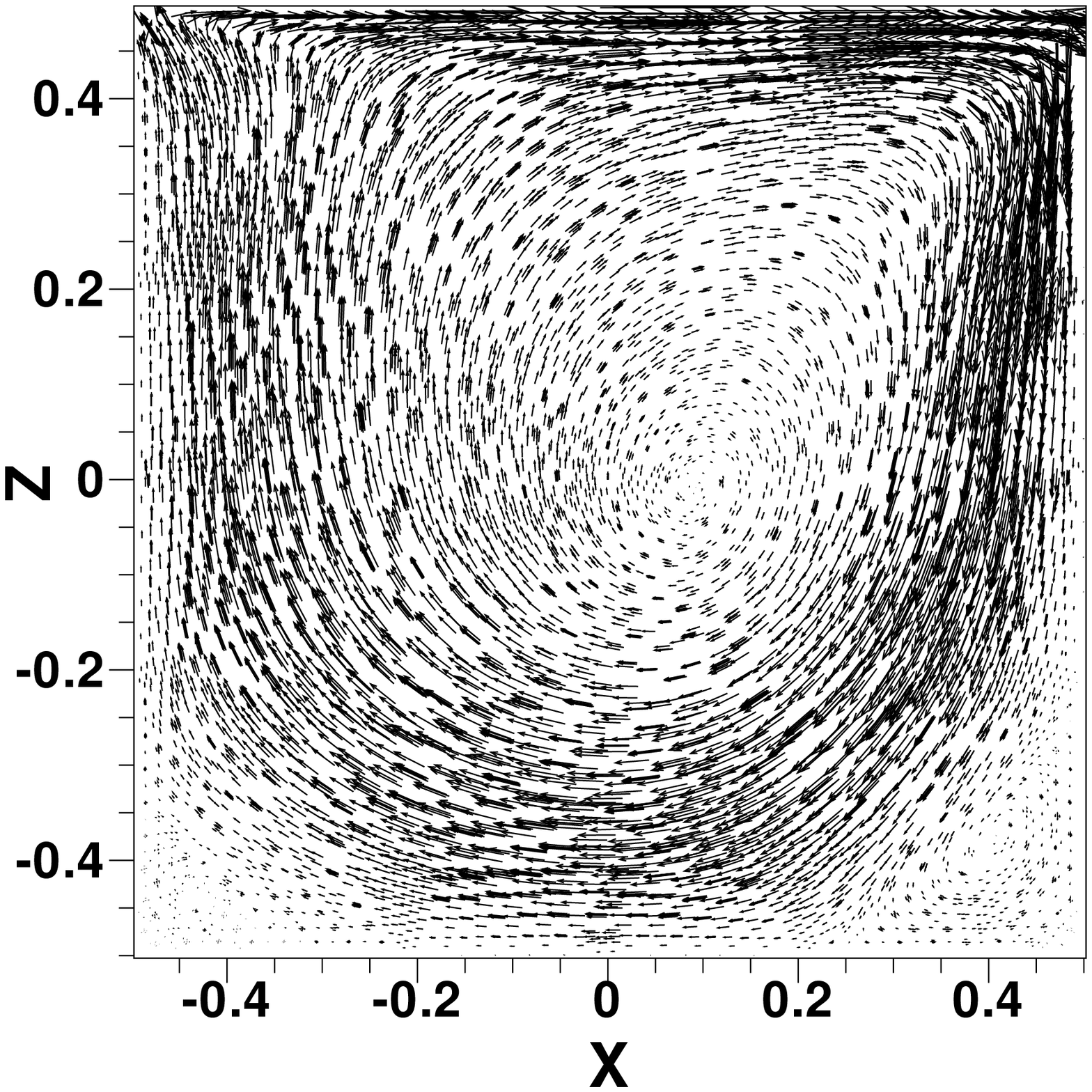}}
  \subfigure{\includegraphics[scale=0.2]{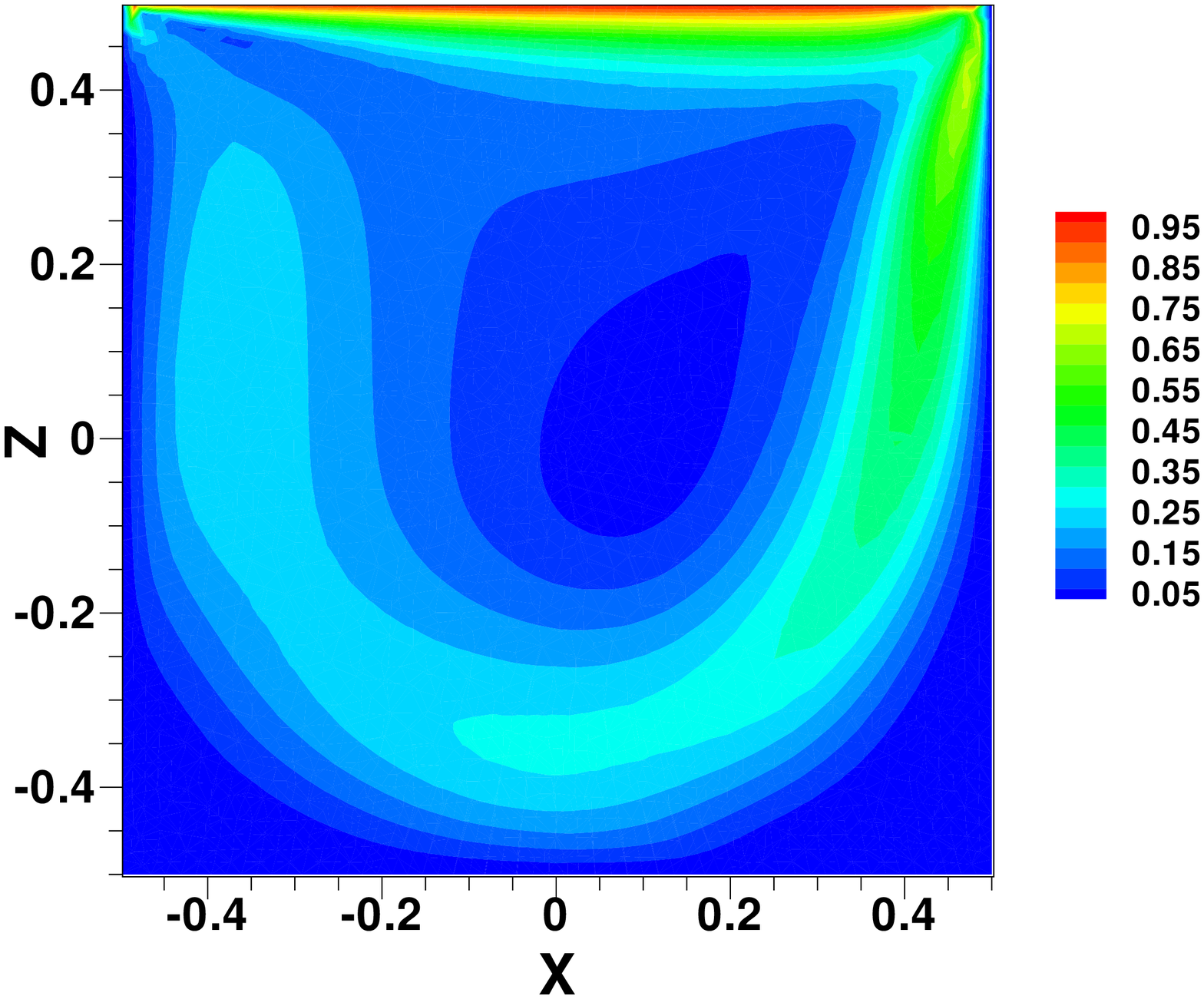}}
  \subfigure{\includegraphics[scale=0.2]{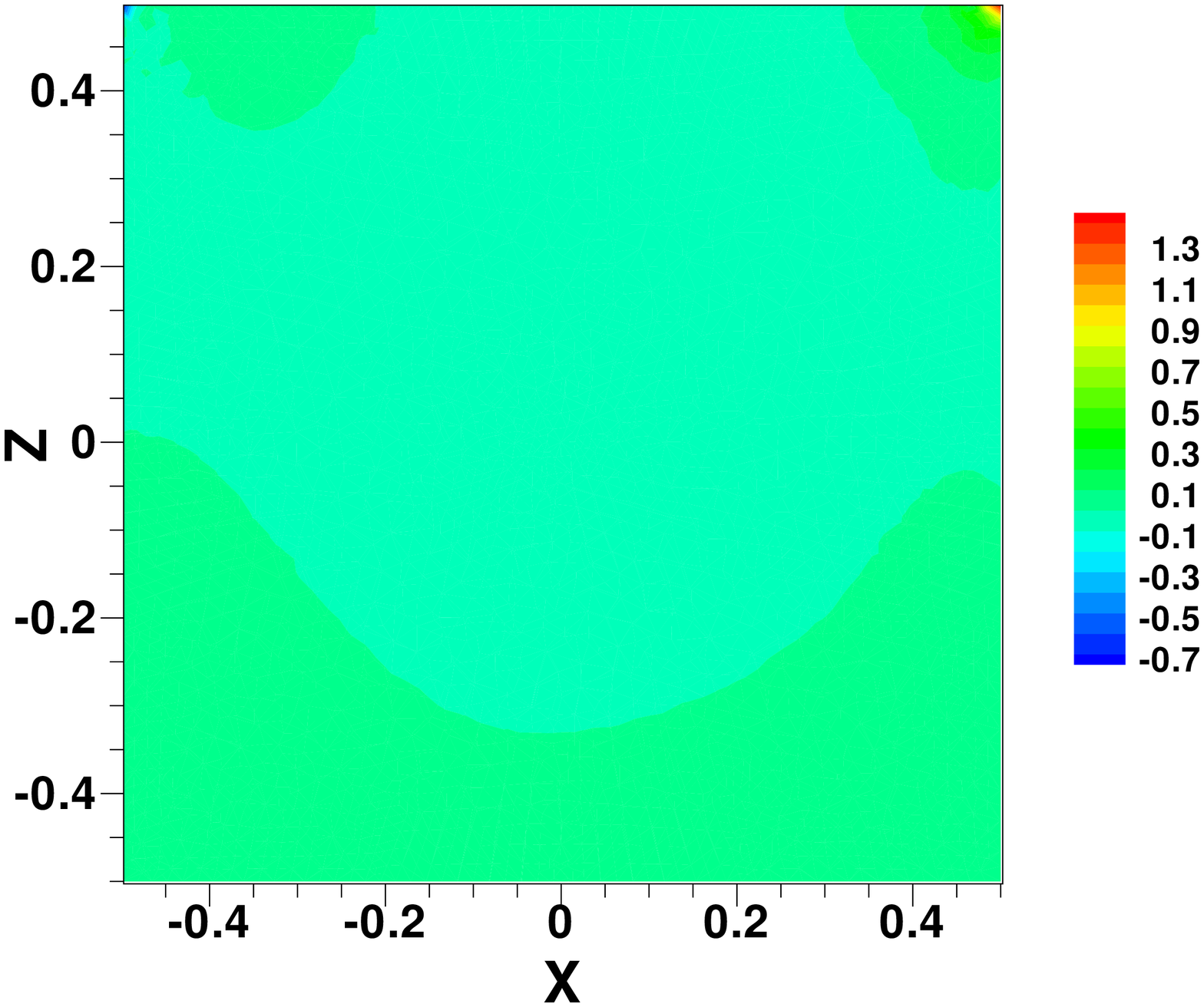}}  
    \subfigure{\includegraphics[scale=0.2]{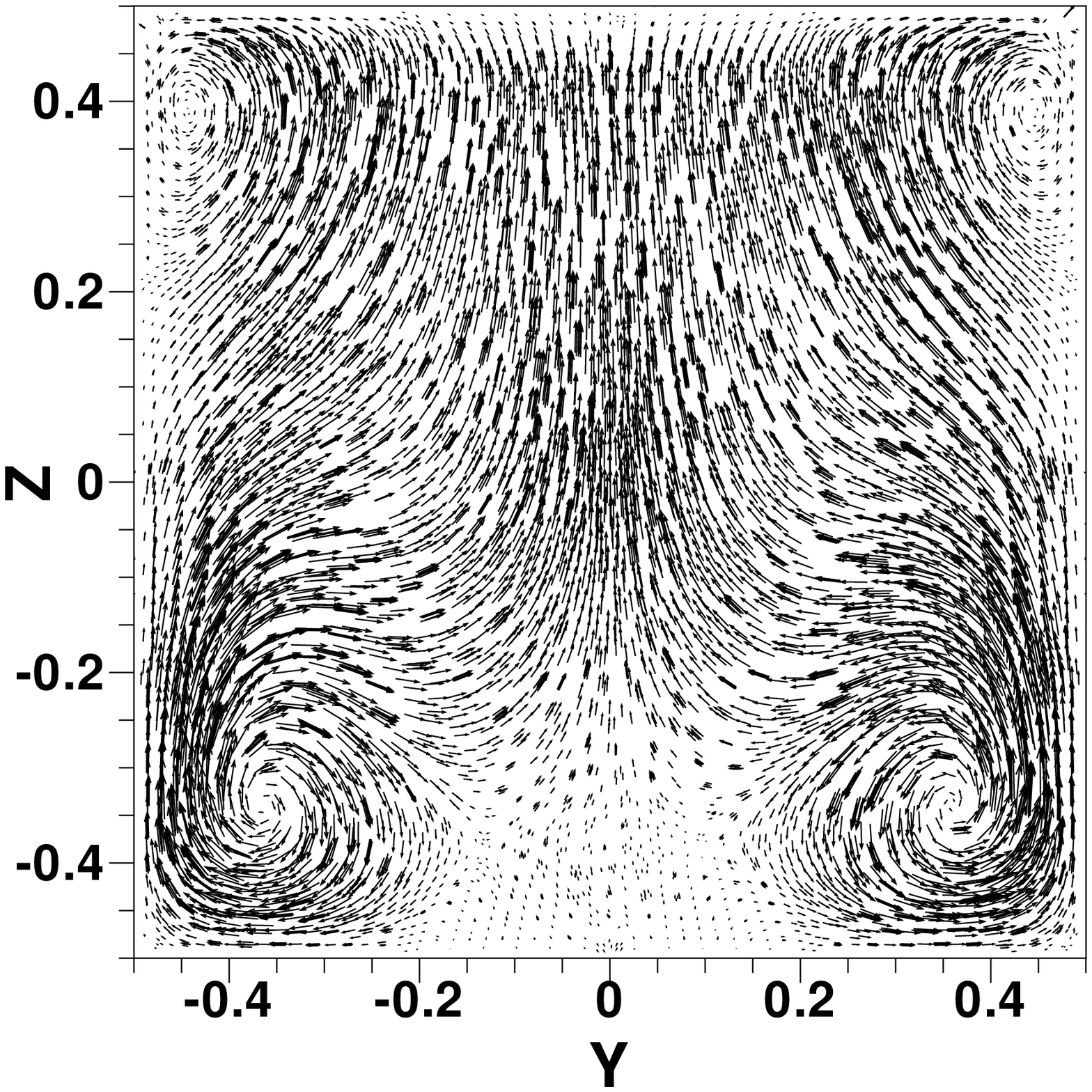}}
  \subfigure{\includegraphics[scale=0.2]{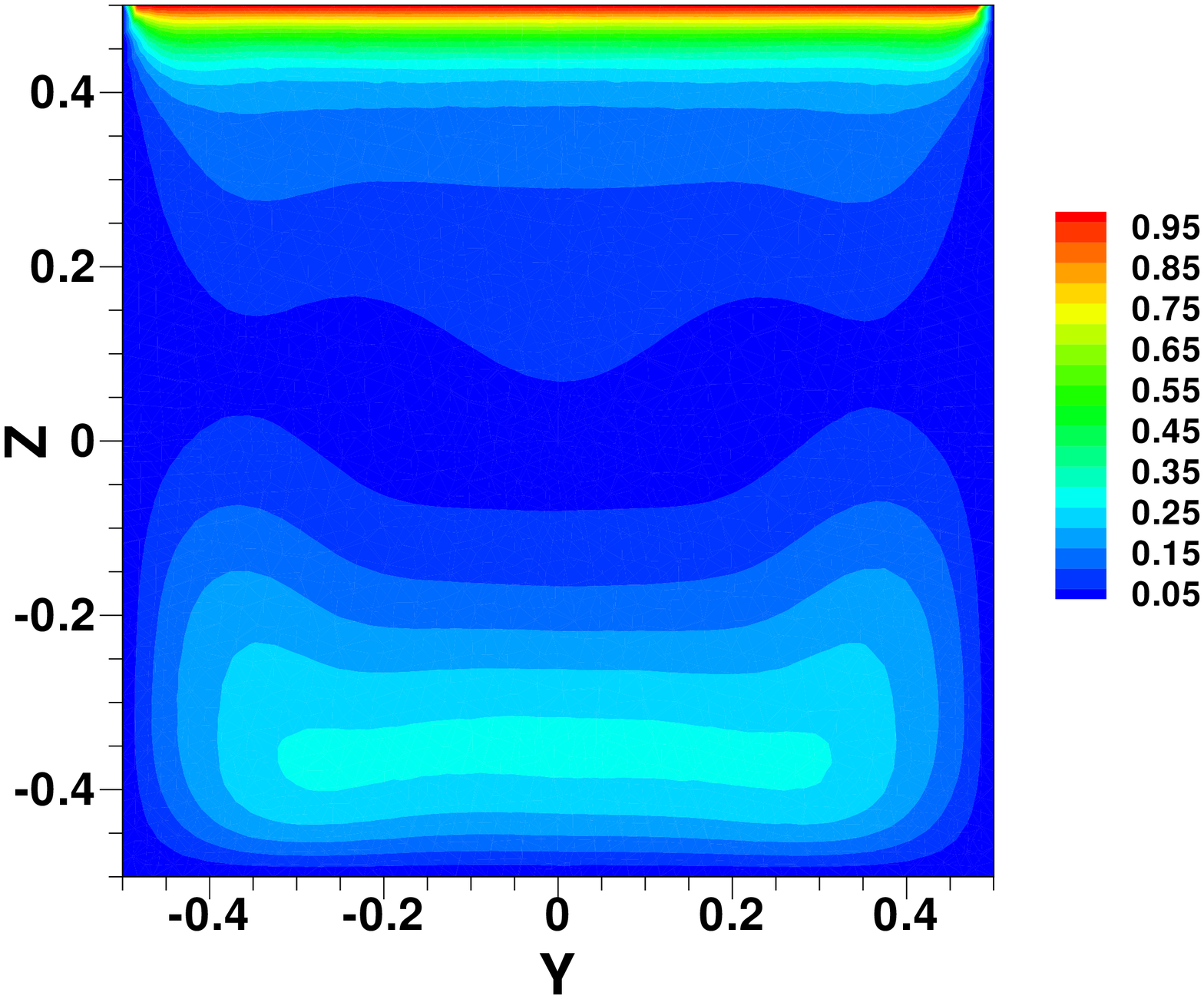}}
  \subfigure{\includegraphics[scale=0.2]{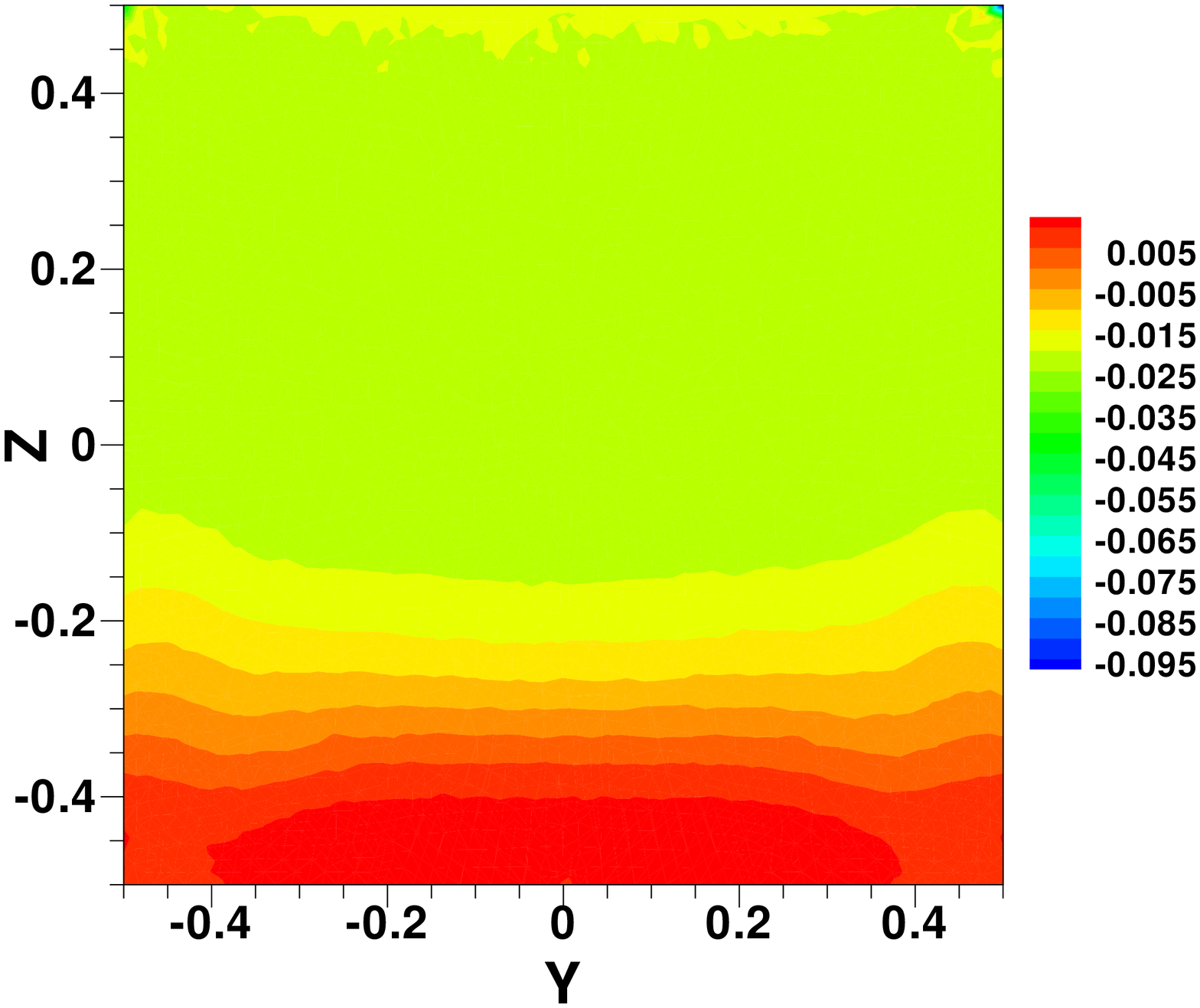}}  
	\caption{3D lid-driven cavity: numerical results. The left column shows vector plots of the velocity, the middle column shows the velocity magnitude, and the right column shows pressure contours in the $x$-$y$ plane (top row), $x$-$z$ plane (middle row), and $y$-$z$ planes (bottom row).}
	\label{fig:CavVelPresVec}
\end{figure}

Figure \ref{fig:CavVelPresVec} shows the results of the three-dimensional lid-driven cavity problem for $Re = 800$ which were obtained by direct solution of the steady state Navier--Stokes equations using the Newton--Schur solution approach.  Notice the presence of the Taylor-G$\ddot{\mathrm{o}}$rtler vortices shown in the bottom corners of the z-y profile. The presence of these vortices at $Re = 800$ may suggest that the transition to turbulent flow may be occurring well before $Re = 2000$. Notice the similarity between the results presented in Figure \ref{fig:CavVelPresVec} and those reported for $Re = 1000$ in \cite{Jiang} which also suggests that unstable flow may be occurring at a lower Reynolds number than previously reported.

\subsection{Three-dimensional jet flow from an orifice}
\begin{figure}[htb!]
	\centering
  \includegraphics[scale=0.5]{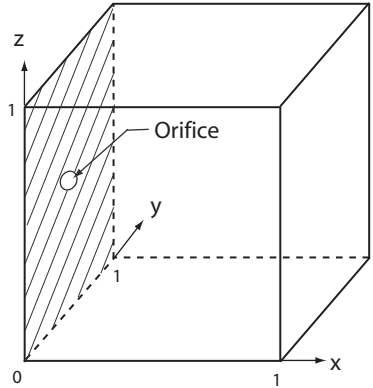}
	\caption{3D jet problem: problem description and boundary conditions. All boundaries are given traction free conditions ($\boldsymbol{\sigma} = 0$) except for the no-slip wall in which the orifice lies. The velocity is prescibed as zero in all directions on the no-slip wall. Over the orifice a parabolic velocity is prescibed in the $x$-direction with maximum magnitude of 1.0 and zero along the edges of the orifice. The velocities in $y$ and $z$ over the orifice are zero.}
	\label{fig:3DJet}
\end{figure}

The three-dimensional jet flow from an orifice problem represents the laminar evolution of vortices formed by flow through a small hole in an infinite no-slip wall. The problem description and boundary conditions are shown in Figure \ref{fig:3DJet}. As the flow is injected into the domain, a mushroom like vortex forms which grows in time. The viscosity was given as 0.001, which corresponds with a Reynolds number of 1000. A time step of 0.01 was also used which is typical for similar problems.

\begin{figure}[htb!]
	\centering
  \includegraphics[scale=0.35]{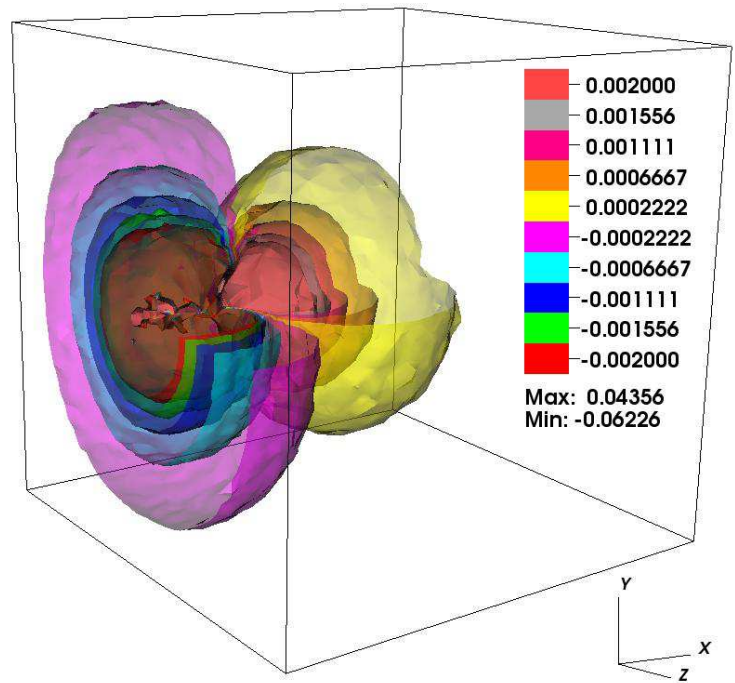}
	\caption{3D jet problem: pressure contours with a slice removed to show detail at time $t = 1.0\mathrm{s}$. The pressure contours are moving forward and spreading sideways. Also note the large region of high pressure that forms in front of a spreading region of low pressure.}
	\label{fig:3DJetPressure}
\end{figure}

\begin{figure}[htb!]
	\centering
  \includegraphics[scale=0.35]{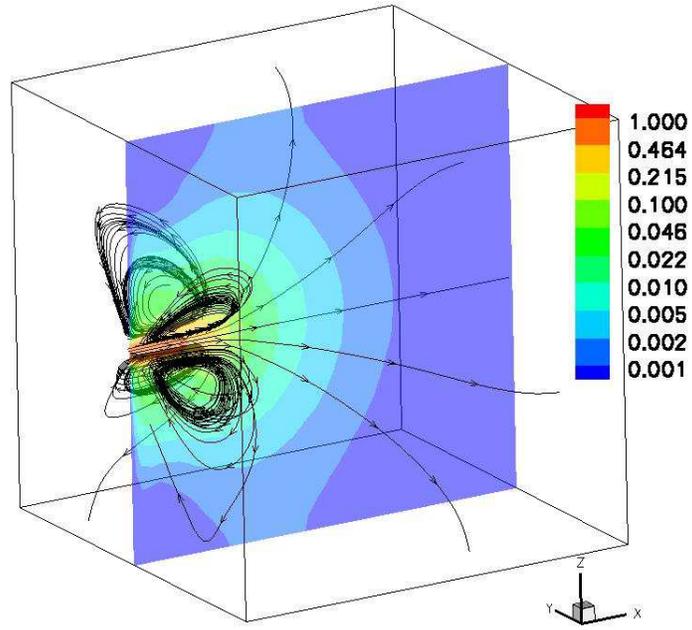}
	\caption{3D jet problem: velocity magnitude contours and streamtraces at time $t = 1.0\mathrm{s}$. The mushroom-shaped vortex formed by the jet at the orifice is clearly shown by the streamtraces.}
	\label{fig:3DJetVelocity}
\end{figure}

Figure \ref{fig:3DJetPressure} shows the pressure contours at time $t = 1.0\mathrm{s}$. A section of the contours have been removed in order to view the inside of the vortex. Figure \ref{fig:3DJetVelocity} shows the velocity magnitude contours and the streamtraces for the three-dimensional jet problem. Notice that the streamtraces are similar to those from the two-dimensional analysis found in \cite{Donea,TurnerNavierStokes}. Notice that the pressure contours are moving forward and spreading out. A large region of high pressure forms in front of the spreading region of low pressure.  The low pressure region corresponds with the vortex shown in the streamtraces of Figure \ref{fig:3DJetVelocity}.

\begin{figure}[htb!]
	\centering
  \subfigure{\includegraphics[scale=0.5]{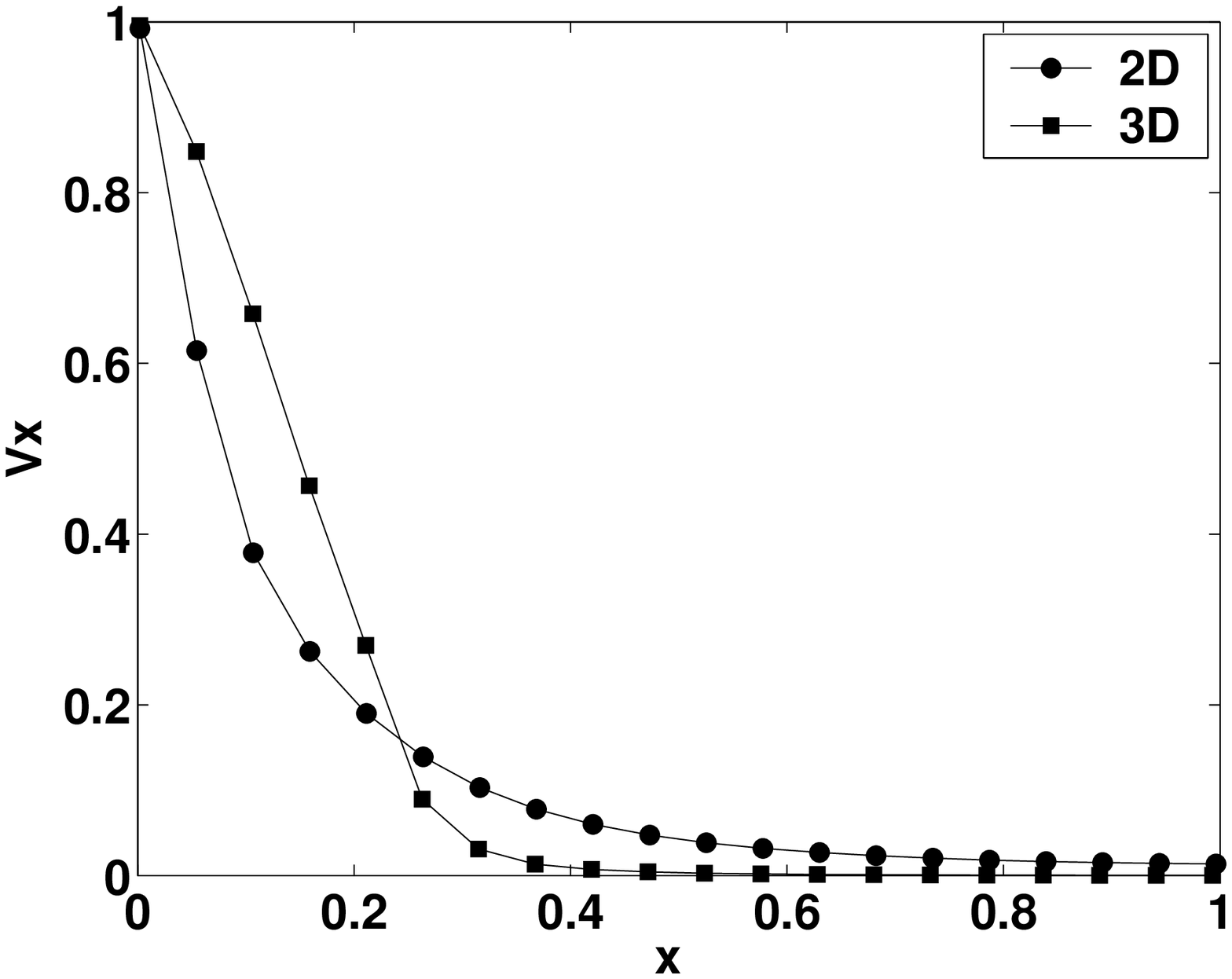}}
  \subfigure{\includegraphics[scale=0.5]{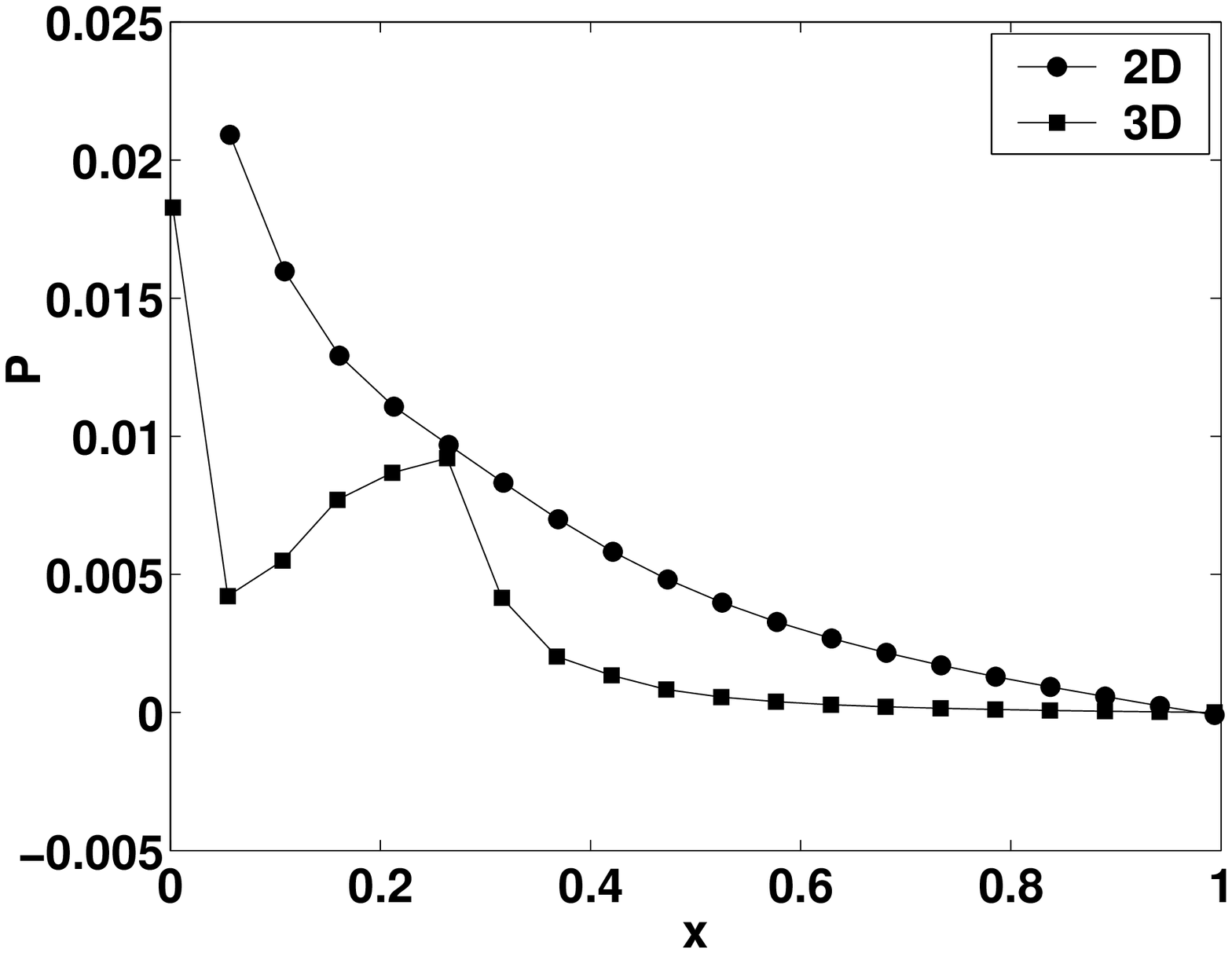}}
	\caption{3D jet problem: a comparison between two-dimensional and three-dimensional results.}
	\label{fig:Jet2D3DCompare}
\end{figure}

In Figure \ref{fig:Jet2D3DCompare} we show a comparison between the two-dimensional results found in \cite{TurnerNavierStokes} and the three-dimensional results presented herein. Notice that for both pressure and velocity the results do not match closely suggesting that three-dimensional effects are present for the jet flow problem. 

\subsection{Parallel speedup and isoefficiency}
To explore the parallel features of the Newton--Schur solution approach a parallel performance study was performed.  The lid-driven cavity problem was solved using a variety of mesh sizes and number of processors.  The boundary conditions for the three-dimensional lid-driven cavity are shown in Figure \ref{fig:3DCavity}. The sizes of each mesh are listed in Table \ref{table:MeshSizes}. The largest mesh consists of over one million degrees of freedom. The problem was solved for each mesh size on up to 128 processors.

\begin{table}[ht!]
\caption{Mesh Sizes in Degrees of Freedom (DOF) for Parallel Results}
\centering
\begin{tabular}{c r r r}
\hline
Mesh  & Coarse DOF & Fine DOF & Total DOF \\
\hline
 A   &  34,591  & 167,289  &   201,880   \\
 B   &  61,697  & 296,496  &   358,193   \\
 C   &  95,803  & 457,476  &   553,279   \\
 D   &  13,5571 & 644,676  &   780,247   \\
 E   &  224,243 & 1,061,157 &   1,285,400  \\
\hline
\end{tabular}
\label{table:MeshSizes}
\end{table}

\begin{figure}[htb!]
	\centering
  \includegraphics[scale=0.5]{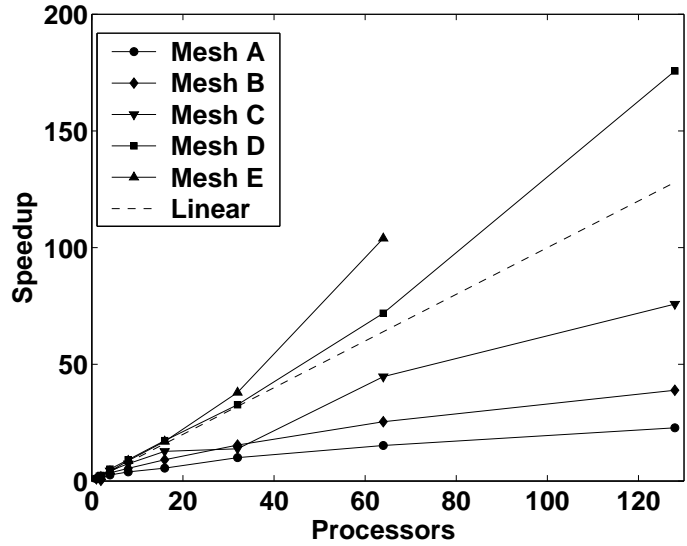}
	\caption{Parallel speedup for various mesh sizes. For the number of degrees of freedom in each mesh, see Table \ref{table:MeshSizes}.}
	\label{fig:ParallelSpeedupNS}
\end{figure}

\begin{figure}[htb!]
	\centering
  \includegraphics[scale=0.5]{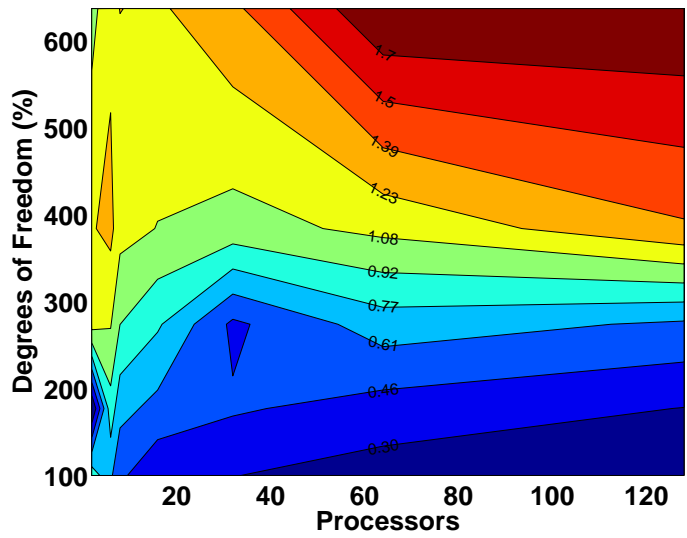}
	\caption{Isoefficiency contours for the Newton-Schur solution approach.}
	\label{fig:IsoE}
\end{figure}

The parallel speedup for each mesh is shown in Figure \ref{fig:ParallelSpeedupNS}. We measure the speedup as the ratio of the parallel execution time over the serial execution time. Notice that Meshes D and E superlinear speedup is obtained, while for Meshes A, B, and C, the speedup is less than linear.  These results correspond with the expected performance as the problem size increases.  For a larger problem size, the ratio of computation time to communication time becomes smaller, allowing for more efficient computation. For a fixed problem size, no algorithm is scalable in the sense that eventually the communication costs will far exceed the computation cost as the number of processors increases. In order to better gage the parallel features of the Newton--Schur solution approach, an isoefficiency study was performed. Isoefficiency contours characterize the problem size required for a given number of processors to achieve a particular value of efficiency. By efficiency, we are referring to the ratio of the serial costs over the parallel costs. For a given serial computation time, $T_1$, and parallel time, $T_p$, for $p$ processors, the efficiency, $E$, is given as
\begin{align}
E_p = \frac{T_1}{pT_p}
\end{align}
Figure \ref{fig:IsoE} shows the isoefficiency contours for each of the meshes in Table \ref{table:MeshSizes}. An algorithm that is poorly scalable will have almost vertical isoefficiency contours, which correspond to a much greater problem size required to obtain a particular efficiency for a given number of processors \cite{HeathNotes}. Notice that for greater than 32 processors, the isoefficiency contours are almost horizontal which suggests that to sustain a certain amount of efficiency for a greater number of processors does not require a substantially larger problem size. The isoefficiency contours in Figure \ref{fig:IsoE} reveal that the Newton--Schur solution approach is very scalable. Good parallel performance is obtained for a reasonable problem size.

\begin{table}[ht!]
\caption{Summary of Parallel Results}
\centering
\begin{tabular}{c c c c}
\hline
     & Max. Efficiency & Theoretical  & Speedup \\
Mesh & Obtained (\%)   & Speedup      & Obtained  \\
\hline
 A   &  86  & 27  &  22  \\
 B   &  86  & 42  &  38  \\
 C   &  92  & 146  &  75  \\
 D   &  122 &  - &  175  \\
 E   &  135 &  - &  103  \\
\hline
\end{tabular}
\label{table:ParallelResults}
\end{table}

The parallel performance data is summarized in Table \ref{table:ParallelResults}.  Values for the theoretical speedup were obtained using Amdahl's law, which estimates the speedup based on the serial fraction of the computation. The speedup cannot exceed the inverse of the serial fraction. For meshes D and E Amdahl's law is not appropriate because the efficiency is greater than 100\%. This may be due to the fact that for such a large problem size, the data does not fit in the cache memory on one processor leading to an extremely long serial computation time. Once the problem has been partitioned over several processors, the data fits into the cache which results in an efficiency greater than 100\%. The superlinear speedup obtained for meshes D and E is also explained by this circumstance.

It is important to note that the number of unknowns is greatly increased by the addition of the fine scale variables.  For example, for the lid-driven cavity problem, the total number of degrees of freedom for the most refined mesh is 1,285,400 whereas the number of coarse-scale degrees of freedom is 224,243 which makes up only 21\% of the total. By using the Schur's complement implementation the initial problem size is reduced by roughly 80\%. The scalability of the algorithm helps to offset the extra computational cost produced by the dual scales. Although the problem size is much greater, the parallel performance is improves as the problem size grows. Also, the fine-scale problem may be further partitioned due to the uncoupled nature of the fine-scale terms.

\begin{remark}
This example also illustrates the increased cost associated with using the variational multiscale framework. The computational cost associated with the variational multiscale framework is approximately four times the cost of solving the coarse-scale problem alone.
\end{remark}

%

\section{CONCLUSIONS}
In summary, we have investigated the performance a Schur's complement implementation of a consistent linearization of the incompressible Navier--Stokes equations for a variety of three-dimensional problems. We obtained high Reynolds solutions to the lid-driven cavity problem without the use of a time-marching scheme. The results obtained for $Re = 800$ suggest that instabilities in the flow, which represent the onset of the transition to turbulent flow may occur at Reynolds numbers less than 1000. We also showed that three dimensional effects may be present for the jet flow through an orifice problem. The results reveal a number of interesting features of the variational multiscale framework. For example, for large problems the dual scales introduce a significant amount of additional degrees of freedom. For the lid-driven cavity problem, the fine scale variable comprised over 80\% of the total degrees of freedom. Regarding the parallel performance of the Newton--Schur solution approach, we were able to achieve superlinear speedup for meshes with greater than one million degrees of freedom. By means of an isoefficiency study, we also showed that the Newton-Schur solution algorithm is scalable for reasonable problem sizes.
%
%

\section*{ACKNOWLEDGMENTS}
The research reported herein was supported by the Computational Science and Engineering Fellowship (D. Z. Turner) and The Department of Energy (K. B. Nakshatrala) through a SciDAC-2 project (Grant No. DOE DE-FC02-07ER64323). This support is gratefully acknowledged. The opinions expressed in this paper are those of the authors and do not necessarily reflect that of the sponsor.

\bibliographystyle{unsrt}
\bibliography{../../DISSERTATION/References/references}

\begin{thebibliography}{10}

\bibitem{Gunzburger}
M.~D. Gunzburger.
\newblock {\em Finite Element Methods for Viscous Incompressible Flows: A Guide
  to Theory Practice, and Algorithms}.
\newblock Academic Press, Inc., San Diego, CA, 1989.

\bibitem{Gresho}
P.~M. Gresho and R.~L. Sani.
\newblock {\em Incompressible Flow and the Finite Element Method Volume 2:
  Isothermal Laminar Flow}.
\newblock John Wiley and Sons, Ltd., New York, 2000.

\bibitem{Donea}
J.~Donea and A.~Huerta.
\newblock {\em Finite Element Methods for Flow Problems}.
\newblock John Wiley and Sons, Ltd., West Sussex, England, 2003.

\bibitem{Hughes2}
T.~J.~R. Hughes.
\newblock Multiscale phenomena: Green's functions, the {D}irichlet-to-{N}eumann
  formulation, subgrid scale models, bubbles and the origins of stabilized
  methods.
\newblock {\em Computer Methods in Applied Mechanics and Engineering},
  127:387--401, 1995.

\bibitem{Jiang}
B.~N. Jiang.
\newblock {\em The Least-Squares Finite Element Method}.
\newblock Springer, New York, 1998.

\bibitem{Russo}
A.~Russo.
\newblock Bubble stabilization of finite element methods for the linearized
  incompressible {N}avier--{S}tokes equations.
\newblock {\em Computer Methods in Applied Mechanics and Engineering},
  132:335--343, 1996.

\bibitem{Masud}
A.~Masud and R.~A. Khurram.
\newblock A multiscale finite element method for the incompressible
  {N}avier-{S}tokes equations.
\newblock {\em Computer Methods in Applied Mechanics and Engineering},
  195:1750--1777, 2006.

\bibitem{Franca2}
L.~P. Franca and A.~Nesliturk.
\newblock On a two-level finite element method for the incompressible
  {N}avier--{S}tokes equations.
\newblock {\em International Journal for Numerical Methods in Engineering},
  52:433--453, 2001.

\bibitem{PSPG}
T.~J.~R. Hughes, L.~P. Franca, and M.~Balestra.
\newblock A new {FEM} for {CFD}: V. {C}ircumventing the {B}abu\v ska-{B}rezzi
  condition: A stable {P}etrov--{G}alerkin formulation of the {S}tokes problem
  accomodating equal order interpolations.
\newblock {\em Computer Methods in Applied Mechanics and Engineering},
  59:85--99, 1986.

\bibitem{SUPG}
T.~J.~R. Hughes and A.~N. Brooks.
\newblock Streamline upwind/{P}etrov--{G}alerkin methods for convection
  dominated flows with emphasis on the incompressible {N}avier-{S}tokes
  equations.
\newblock {\em Computer Methods in Applied Mechanics and Engineering},
  32:199--259, 1982.

\bibitem{TurnerStokes}
D.~Z. Turner, K.~B. Nakshatrala, and K.~D. Hjelmstad.
\newblock On the stability of bubble functions and a mixed finite element
  formulation for the {S}tokes problem.
\newblock {\em International Journal for Numerical Methods in Fluids},
  Submitted, 2008.

\bibitem{TurnerNavierStokes}
D.~Z. Turner, K.~B. Nakshatrala, and K.~D. Hjelmstad.
\newblock Consistent {N}ewton-{R}aphson vs. fixed-point for variational
  multiscale formulations for incompressible {N}avier--{S}tokes.
\newblock {\em International Journal for Numerical Methods in Fluids},
  Submitted, 2008.

\bibitem{Nakshatrala}
K.~B. Nakshatrala, A.~Masud, and K.~D. Hjelmstad.
\newblock On finite element formulations for nearly incompressible linear
  elasticity.
\newblock {\em Computational Mechanics}, 41:547--561, 2008.

\bibitem{Nakshatrala2}
K.~B. Nakshatrala, D.~Z. Turner, K.~D. Hjelmstad, and A.~Masud.
\newblock A stabilized mixed finite element method for {D}arcy flow based on a
  multiscale decomposition of the solution.
\newblock {\em Computer Methods in Applied Mechanics and Engineering},
  195(33-36):4036--4049, 2006.

\bibitem{Masud2}
A.~Masud and R.~A. Khurram.
\newblock A multiscale/stabilized finite element method for the
  advection-diffusion equation.
\newblock {\em Computer Methods in Applied Mechanics and Engineering},
  193:1997--2018, 2004.

\bibitem{Koobus}
B.~Koobus and C.~Farhat.
\newblock A variational multiscale method for the large eddy simulation of
  compressible turbulent flows on unstructured meshes--applications to vortex
  shedding.
\newblock {\em Computer Methods in Applied Mechanics and Engineering},
  193:1367--1383, 2008.

\bibitem{Hughes3}
T.~J.~R. Hughes, L.~Mazzei, and K.~E. Jansen.
\newblock Large eddy simulation and the variational multiscale method.
\newblock {\em Computing and Visualization in Science}, 3:47--59, 2000.

\bibitem{Kulkarni}
D.~V. Kulkarni and D.~A. Tortorelli.
\newblock {\em A domain decomposition based two-level {N}ewton scheme for
  nonlinear problems.}
\newblock Domain Decomposition Methods in Science and Engineering. Springer,
  2004.

\bibitem{Kulkarni2}
D.~V. Kulkarni, D.~A. Tortorelli, and M.~Wallin.
\newblock A newton--schur alternative to the consistent tangent approach in
  computational plasticity.
\newblock {\em Computer Methods in Applied Mechanics and Engineering},
  196:1169--1177, 2007.

\bibitem{Brezzi}
F.~Brezzi and M.~Fortin.
\newblock {\em Mixed and Hybrid Finite Element Methods}, volume 15 of Springer
  Series in Computational Mathematics.
\newblock Springer-Verlag, New York, USA, 1997.

\bibitem{Arnold}
D.~N. Arnold, F.~Brezzi, and M.~Fortin.
\newblock A stable finite element for the {S}tokes equations.
\newblock {\em Estratto da Calcolo}, 21(4):337--344, 1984.

\bibitem{Graham}
A.~Graham.
\newblock {\em Kronecker Products and Matrix Calculus: With Applications}.
\newblock Halsted Press, 1981.

\bibitem{Turing}
University of~Illinois~at Urbana-Champaign: Computational~Science and
  Engineering.
\newblock {Turing Cluster} {W}eb page, 2008.
\newblock http://www.cse.uiuc.edu/turing/.

\bibitem{petsc-web-page}
S.~Balay, K.~Buschelman, W.D. Gropp, D.~Kaushik, M.G. Knepley, L.C. McInnes,
  B.F. Smith, and H.~Zhang.
\newblock {PETSc} {W}eb page, 2001.
\newblock http://www.mcs.anl.gov/petsc.

\bibitem{Tecplot}
Inc. Tecplot.
\newblock Tecplot 360, 2008.
\newblock http://www.tecplot.com.

\bibitem{Visit}
Lawrence Livermore~National Laboratory.
\newblock Vis{I}t 1.8.1, 2005.
\newblock http://www.llnl.gov/visit.

\bibitem{Koseff}
J.~R. Koseff and R.~L. Street.
\newblock On end wall effects in a lid-driven cavity flow.
\newblock {\em Journal of Fluids Engineering}, 106:385--389, 1984.

\bibitem{Saric}
W.~S. Saric.
\newblock G$\ddot{\mathrm{o}}$rtler vortices.
\newblock {\em Annual Review of Fluid Mechanics}, 26:379--409, 1994.

\bibitem{HeathNotes}
M.~Heath.
\newblock Chapter 4: Parallel performance, 2006.
\newblock University of Illinois: CS 554 Course Notes.

\end{thebibliography}
\newpage
%

\end{document}